\input amstex\documentstyle {amsppt}  
\pagewidth{12.5 cm}\pageheight{19 cm}\magnification\magstep1
\topmatter
\title Cuspidal local systems and graded Hecke algebras, III\endtitle
\author G. Lusztig\endauthor
\address Department of Mathematics, M.I.T., Cambridge, MA 02139\endaddress
\thanks{Supported in part by the National Science Foundation}\endthanks
\abstract{We prove a strong induction theorem and classify the tempered and
square integrable representations of graded Hecke algebras.}\endabstract
\endtopmatter
\document     
\define\lan{\langle}
\define\ran{\rangle}
\define\bsl{\backslash}
\define\mto{\mapsto}
\define\lra{\leftrightarrow}

\define\sm{\smallmatrix}
\define\esm{\endsmallmatrix}
\define\sub{\subset}
\define\bxt{\boxtimes}
\define\tim{\times}
\define\ti{\tilde}
\define\nl{\newline}
\redefine\i{^{-1}}
\define\fra{\frac}
\define\emp{\emptyset}
\define\frl{\forall}
\define\un{\underline}

\define\ot{\otimes}
\define\bbq{\bar{\bq}_l}

\define\ad{\text{\rm ad}}
\define\Ad{\text{\rm Ad}}
\define\Hom{\text{\rm Hom}}

\define\Ker{\text{\rm Ker}}
\redefine\Im{\text{\rm Im}}

\define\tr{\text{\rm tr}}

\define\opl{\oplus}
\define\sha{\sharp}

\define\cok{\text{\rm coker}}

\define\al{\alpha}

\define\ga{\gamma}
\define\de{\delta}
\define\ep{\epsilon}
\define\io{\iota}

\define\rh{\rho}
\define\si{\sigma}
\define\ta{\tau}

\define\la{\lambda}
\define\ze{\zeta}
\define\ph{\phi}
\define\ps{\psi}
\define\vp{\varpi}

\define\Ga{\Gamma}
\define\De{\Delta}

\define\La{\Lambda}

\define\boc{\bold c}

\define\boo{\bold o}

\define\bor{\bold r}
\define\btt{\bold t}

\define\bc{\bold C}

\define\bh{\bold H}

\define\bn{\bold N}

\define\bq{\bold Q}
\define\br{\bold R}
\define\bs{\bold S}
\define\bt{\bold T}

\define\bz{\bold Z}

\define\ca{\Cal A}
\define\cb{\Cal B}
\define\cc{\Cal C}

\define\ce{\Cal E}

\define\cg{\Cal G}

\define\ci{\Cal I}
\define\cj{\Cal J}
\define\ck{\Cal K}
\define\cl{\Cal L}
\define\cm{\Cal M}
\define\cn{\Cal N}
\define\co{\Cal O}
\define\cp{\Cal P}

\define\car{\Cal R}
\define\cs{\Cal S}
\define\ct{\Cal T}

\define\cx{\Cal X}
\define\cy{\Cal Y}

\define\fa{\frak a}
\define\fb{\frak b}
\define\fc{\frak c}

\define\fg{\frak g}
\define\fh{\frak h}

\define\fl{\frak l}

\define\fn{\frak n}

\define\fp{\frak p}
\define\fq{\frak q}

\define\fs{\frak s}
\define\ft{\frak t}
\define\fu{\frak u}

\define\fz{\frak z}

\define\pa{\frak P}
\define\bap{\bar P}
\define\ubap{\un{\bap}}
\define\uzp{\un{Z_{\bap}}}
\define\baq{\bar Q}
\define\bfu{\bar\fu}
\define\bam{\bar M}
\define\bau{\bar U}
\define\bae{\bar E}
\define\bacg{\bar{\cg}}

\define\baD{\bar\De}
\define\bad{\bar\de}
\define\bct{\bar\ct}
\define\bax{\bar x}

\define\tbu{\ti{\bar u}}
\define\tbx{\ti{\bar x}}

\define\ug{\un G}
\define\up{\un P}
\define\uup{\un{U_P}}
\define\uuq{\un{U_Q}}
\define\uq{\un Q}
\define\ut{\un T}
\define\ul{\un L}
\define\uz{\un Z}
\define\uxi{\un\xi}
\define\upi{\un\pi}
\define\ubar{\un{\bar R}}
\define\dfg{\dot{\fg}}
\define\fzg{\fz_\fg}

\define\tm{\ti m}
\define\tx{\ti X}
\define\ty{\ti y}
\define\tf{\ti f}
\define\do{\dot{}}
\define\dx{\dot X}
\define\dy{\dot Y}
\define\dcl{\dot{\cl}}
\define\dq{\dot\fq}
\define\da{\dot A}
\define\dca{\dot{\ca}}
\define\impl{\implies}
\define\Irr{\text{\rm Irr}}
\define\et{\eta}
\define\bocp{\boc_P}
\define\BW{BW}
\define\KL{KL}
\define\LE{Le}
\define\LSQ{L1}
\define\LIC{L2} 
\define\LCS{L3} 
\define\LI{L4}
\define\LII{L5}
\define\LPE{L6}
\define\LUN{L7}  
\define\LB{L8}
\define\RE{R}  
\define\WA{W}

\head Introduction\endhead
Let $\cg$ be the group of rational points of a simple adjoint algebraic group 
over a $p$-adic field, which is an inner form of a split group. Consider the 
set of isomorphism classes of irreducible admissible representations of $\cg$ 
whose restriction to some parahoric subgroup contains some irreducible
unipotent cuspidal representation of that parahoric subgroup modulo its 
"unipotent radical". The classification of such "unipotent" representations of
$\cg$ has been established in \cite{\LUN} (see also \cite{\KL} for an earlier 
special case) in accordance with a conjecture of Langlands (refined in
\cite{\LSQ}). In the special case considered in \cite{\KL}, the tempered and
square-integrable representations were also explicitly described; the main 
tool to
do so was an "induction theorem" \cite{\KL, 6.2} for affine Hecke algebras 
with equal parameters. But in the context of \cite{\LUN} the induction theorem
was missing
and the problem of describing explicitly the unipotent representations that are
tempered or square integrable representations remained open.

One of the techniques used in \cite{\LUN} was the reduction of 
the (equivalent) classification problem for certain affine Hecke algebras with
unequal parameters to the problem of classifying the simple modules of certain
"graded" Hecke algebras which could be done using methods of equivariant 
homology. By these methods one can reduce the problem of describing the
tempered or square integrable unipotent representations to the
analogous problem for graded Hecke algebras. This last problem is solved in the
present paper. As in \cite{\KL} one of the key ingredients is an "induction
theorem". In fact we will prove a strong form of the induction theorem 
(without "denominators") inspired by \cite{\LB, 7.11} which implies 
the classification of tempered and square integrable representations.

In the case where $\cg$ has small rank, the classification of square integrable
unipotent representations of $\cg$ has been given in \cite{\RE}.

After this work was completed, I received the preprint \cite{\WA} where most 
results of the present paper (but not the strong induction theorem) are 
obtained independently in the case $\cg=SO_{2n+1}$.

\head Contents\endhead
1. Preliminaries and statement of results.

2. A strong induction theorem and a proof of Theorem 1.17.

3. Proof of Theorems 1.15, 1.21, 1.22.

\head 1. Preliminaries and statement of results\endhead
\subhead 1.1\endsubhead
Unless otherwise specified, all algebraic varieties are assumed to be over
$\bc$. If $X$ is a subvariety of $X'$ we write $cl(X)$ for the closure of $X$
in $X'$. For a Lie algebra $\fg$ let $\fzg$ be the centre of $\fg$. If $A$ is a
subset of $\fg$, we set $\fz(A)=\{x\in\fg;[x,y]=0\quad\frl y\in A\}$; if 
$\fa$ is a
Lie subalgebra of $\fg$, we set $\fz_\fa(A)=\fz(A)\cap\fa$.

For any algebraic group $\cg$ let $\cg^0$ be the identity component of $\cg$, 
let $U_\cg$ be the unipotent radical of $\cg^0$ and let $Z_\cg$ be the centre 
of $\cg$. Let $\bacg=\cg/U_{\cg}$ and let $\pi_\cg:\cg@>>>\bacg$ be the
canonical homomorphism. Let $\un{\cg}$ be the Lie algebra of $\cg$. If $A$ is a
subset of $\un{\cg}$, we set $Z(A)=\{g\in G;\Ad(g)y=y\quad\frl y\in A\}$. Let
$\upi_\cg:\un{\cg}@>>>\un{\bacg}$ be the map induced by $\pi_\cg$. Let
$\cg_{der}$ be the derived group of $\cg$ (a closed subgroup if $\cg$ is 
connected). If $x\in\un{\cg}$ is a semisimple element, we denote by 
$\lan x\ran$ the smallest torus in $\cg$ whose Lie
algebra contains $x$.

\subhead 1.2\endsubhead
Let $G$ be a connected reductive algebraic group. Let $\fg=\ug$; let $\fg_N$ be
the variety of nilpotent elements of $\fg$. 
Let $\fg_{ss}$ be the set of semisimple elements of $\fg$.
Let $\pa$ be the variety of 
parabolic subgroups of $G$.

A {\it cuspidal datum} for $G$ is a triple $(\cp,\boc,\cl)$ where $\cp$ is a 
$G$-orbit on $\pa$, $\boc$ is a $G$-orbit on the set of pairs $(x,P)$ with 
$P\in\cp$, $x\in\ubap$ is nilpotent, and $\cl$ is an irreducible 
$G$-equivariant local system on $\boc$ such that for some (or any) $P\in\cp$, 
the restriction of $\cl$ to the $\bap$-orbit 
$$\bocp=\{x\in\ubap;(x,P)\in\boc\}$$
(a local system that is automatically $\bap$-equivariant and irreducible) is
cuspidal in the sense of \cite{\LI, 2.2}. 

A {\it cuspidal triple} in $G$ is a triple $(L,C,\ce)$ where $L$ is a Levi 
subgroup of a parabolic subgroup of $G$, $C$ is a nilpotent $L$-orbit in $\ul$
and $\ce$ is an irreducible $L$-equivariant local system on $C$ which is
cuspidal in the sense of \cite{\LI, 2.2}.

To a cuspidal datum $(\cp,\boc,\cl)$ we attach a cuspidal triple as follows: 
let $P\in\cp$, let $L$ be a Levi subgroup of $P$, let $C$ be the nilpotent 
orbit in $\ul$ corresponding to $\bocp$ under the obvious isomorphism 
$\ul@>\sim>>\ubap$ and let $\ce$ be the local system on $C$ corresponding to
$\cl|_{\bocp}$ under the obvious isomorphism $C@>\sim>>\bocp$. Then $(L,C,\ce)$
is a cuspidal triple in $G$. Using \cite{\LII, 6.8(b),(c)}, we see that, 
conversely, any cuspidal triple is obtained as above from a cuspidal datum 
$(\cp,\boc,\cl)$ where $\cp,\boc$ are unique and $\cl$ is unique up to 
isomorphism.

\subhead 1.3\endsubhead
Let $(\cp,\boc,\cl)$ be a cuspidal datum for $G$ and let $Q\in\pa$ be such that
$Q$ contains some $P\in\cp$. Then there is an induced cuspidal datum 
$(\cp',\boc',\cl)$ for $\baq$, defined as follows. Let $\cp'$ be the set of all
subgroups $P'$ of $\baq$ such that $\pi_Q\i(P')\in\cp$. Let $\boc'$ be the set
of all pairs $(x',P')$ where $P'\in\cp'$ and 
$x'\in\ubap'=\un{\overline{\pi_Q\i(P')}}$ is such that $(x',\pi_Q\i(P'))\in\boc$. 
The inverse image of $\cl$ under the map $\boc'@>>>\boc$ given by
$(x',P')\mto(x',\pi_Q\i(P'))$ is denoted again by $\cl$. Then 
$(\cp',\boc',\cl)$ is a cuspidal datum for $\baq$. 

\subhead 1.4\endsubhead
{\it In the remainder of this paper we fix a cuspidal datum $(\cp,\boc,\cl)$ 
for $G$.} 

For each $P\in\cp$ we form the torus $P/P_{der}$ (resp. the vector space 
$\up/[\up,\up]$). If $P,P'\in\cp$, there is a canonical isomorphism 
$P'/P'_{der}@>\sim>>P/P_{der}$ 
(resp. \linebreak $\up'/[\up',\up']@>\sim>>\up/[\up,\up]$) 
induced by $\Ad(g)$ where $g\in G$ is such that $\Ad(g)P'=P$. This is 
independent of the choice of $g$. Hence we may identify $P/P_{der}$ (resp. 
$\up/[\up,\up]$) for any $P\in\cp$ with a single torus $\bt$ (resp. a single 
$\bc$-vector space $\fh$). Thus, for any $P\in\cp$ we have a canonical 
isomorphism $P/P_{der}@>\sim>>\bt$ (resp. $\up/[\up,\up]@>\sim>>\fh$). Since 
for $P\in\cp$ we have $\up/[\up,\up]=\un{P/P_{der}}$, we have canonically 
$\fh=\un{\bt}$.

\subhead 1.5\endsubhead
The set 

$\{P'\in\pa;P' \text{ contains strictly some $P\in\cp$ and is minimal 
with this property}\}$ 
\nl
decomposes into $G$-orbits $(\cp_i)_{i\in I}$. Here $I$ is a finite indexing 
set.

For any $J\sub I$ let $\cp_J$ be the set of all $P'\in\pa$ such that $P'$ 
contains some member of $\cp$ and, for $i\in I$, $P'$ contains some member of 
$\cp_i$ if and only if $i\in J$. Then $\cp_J$ is a $G$-orbit on $\pa$. We have
$\cp=\cp_\emp$, $\cp_i=\cp_{\{i\}}$ for $i\in I$. 

The diagonal action of $G$ on $\cp\tim\cp$ has only finitely many orbits; an
orbit is said to be {\it good}
 if it consists of pairs $(P,P')$ such that $P,P'$ have
a common Levi subgroup. Let $W$ be the set of good $G$-orbits on $\cp\tim\cp$.
There is a natural group structure on $W$ (see \cite{\LII, 7.3}).

\subhead 1.6\endsubhead
Let $P\in\cp$ and let $L$ be a Levi subgroup of $P$. For any $J\sub I$ let 
$P_J$ be the unique member of $\cp_J$ that contains $P$. For $i\in I$ write 
$P_i$ instead of $P_{\{i\}}$. We have $P_{\emp}=P$. Let $T=Z_L^0$; then 
$\ut=\fz_{\ul}$. Let $N(T)$ be the normalizer of $T$ in $G$. Then $W(L)=N(T)/L$
acts naturally (and faithfully) on $\ut$ and on $\ut^*$. We have 
$\fg=\opl_{\al\in\ut^*}\fg^\al$ where 
$$\fg^\al=\{x\in\fg;[y,x]=\al(y)x \quad\frl y\in\ut\}.$$ 
Note that $\fg^\al$ is an $\ul$-module by the $\ad$ action. Let 
$R=\{\al\in\ut^*;\al\ne 0,\fg^\al\ne 0\}$. 

By \cite{\LI, 2.5}, $R$ is a (not necessarily reduced) root system in $\ut^*$ 
with Weyl group $W(L)$. (We do not have to specify the set of coroots since 
they are determined uniquely by $R$ and the Weyl group action on $\ut^*$.) Let
$L_i$ be the Levi subgroup of $P_i$ that contains $L$. There is a unique 
$\al_i\in R$ such that $\fg^{\al_i}\sub\uup$ and 
$\un{L_i}=\opl_{n\in\bz}\fg^{n\al_i}$. Then $\{\al_i;i\in I\}$ is a set of 
simple roots for $R$.

Let $\cc$ be the nilpotent $L$-orbit in $\ul$ which corresponds to $\bocp$ 
under the obvious isomorphism $\ul@>\sim>>\ubap$. Let $y\in\cc$. For $i\in I$
let $c_i$ be the integer $\ge 2$ such that

$\ad(y)^{c_i-2}:\fg^{\al_i}\opl\fg^{2\al_i}@>>>\fg^{\al_i}\opl\fg^{2\al_i}$ is 
$\ne 0$,

$\ad(y)^{c_i-1}:\fg^{\al_i}\opl\fg^{2\al_i}@>>>\fg^{\al_i}\opl\fg^{2\al_i}$ is 
$0$.
\nl
Then $c_i$ is independent of the choice of $P,L$ and $y$.

We identify $W(L)$ with $W$ by $n\mto G-\text{orbit of } (P,nPn\i)$. Via the
obvious isomorphism $\ut@>\sim>>\up/[\up,\up]=\fh$, the $W(L)$ action on $\ut$ 
and $\ut^*$ becomes a $W$ action on $\fh$ and $\fh^*$ (independent of the 
choice of $P,L$) and the vectors $\al_i (i\in I)$ in $\ut^*$ become vectors in
$\fh^*$, denoted again by $\al_i$ (these are also independent of the choice of
$P,L$). The action of $W$ on $\fh^*$ is denoted by $w,\xi\mto{}^w\xi$. For 
$i\in I$ let $s_i$ be the unique element of $W$ which is a reflection in 
$\fh^*$ such that $s_i(\al_i)=-\al_i$. Then $W$ together with $s_i (i\in I)$ is
a Coxeter group. For $J\sub I$ let $W_J$ be the subgroup of $W$ generated by 
$\{s_i;i\in J\}$.

For future use we note the following property:

(a) {\it Let $i\in I$ and let $Q\in\pa$ be such that $P\sub Q$,
$\un{P_i}\not\sub\uq$. Then $\un{L_i}\cap\uq=\opl_{n\in\bn}\fg^{n\al_i}$. In 
particular}, $(\fg^{-\al_i}\opl\fg^{-2\al_i})\cap\uq=0$.

\subhead 1.7\endsubhead
Let $\bh$ be the associative $\bc$-algebra defined by the generators $\uxi$ (in
1-1 correspondence with the elements $\xi\in\fh^*$), $s_i$ (indexed by 
$i\in I$) and $\bor$, subject to the following relations:

(a) $\un{a\xi+a'\xi'}=a\uxi+a'\uxi'$ for any $\xi,\xi'\in\fh^*$ and any
$a,a'\in\bc$;

(b) $\uxi\uxi'=\uxi'\uxi$ for any $\xi,\xi'\in\fh^*$;

(c) $s_i (i\in I)$ satisfy the relations of $W$;

(d) $s_i\uxi-\un{{}^{s_i}\xi}s_i=c_i\fra{\xi-{}^{s_i}\xi}{\al_i}\bor$ for any
$\xi\in\fh^*$ and any $i\in I$.

(e) $\bor$ is central.
\nl
(In (d) we have $\fra{\xi-{}^{s_i}\xi}{\al_i}\in\bc$.) This is the same as the
algebra denoted by $\bh$ in \cite{\LI, 6.3}.

\subhead 1.8\endsubhead
Let 
$$\dfg=\{(y,P)\in\fg\tim\cp;y\in\up;\upi_P(y)\in\bocp+\uzp\}.$$
Let $\pi:\dfg@>>>\fg$ be the first projection. Now $G\tim\bc^*$ acts on $\fg$ 
by $(g,\la):y\mto\la^{-2}\Ad(g)y$, on $\cp$ by $(g,\la):P\mto gPg\i$ and on 
$\dfg$ by $(g,\la):(y,P)\mto(\la^{-2}\Ad(g)y,gPg\i)$. For $y\in\fg_N$ we denote
by $M(y)$ or $M_G(y)$ the stabilizer of $y$ in $G\tim\bc^*$. Thus, 
$$M(y)=\{(g,\la)\in G\tim\bc^*;\Ad(g)y=\la^2y\}.$$
We also have an action of $G\tim\bc^*$ on $\boc$ given by 
$(g,\la):(x,P)\mto(\la^{-2}\Ad(g)x,gPg\i)$. Here we regard $\Ad(g)$ as a map 
$\bocp@>>>\boc_{gPg\i}$. The local system $\cl$ on $\boc$ is automatically
$G\tim\bc^*$-equivariant \cite{\LII, 7.15}. Let $s:\dfg@>>>\boc$ be given by
$s(y,P)=(y',P)$ where $y'\in\bocp,\upi_P(y)-y'\in\uzp$. Then 
$\dcl=s^*\cl$ is a $G\tim\bc^*$-equivariant local system on $\dfg$ and 
$\ck=\pi_!(\dcl^*)$ is (up to shift) a $G\tim\bc^*$-equivariant perverse sheaf
on $\fg$, with a canonical action of $W$, \cite{\LI, 3.4}.

\subhead 1.9\endsubhead
Let $X$ be an algebraic variety with a given morphism $X@>m>>\fg$. Define 
$X@<m'<<\dx@>\dot m>>\dfg$ by the cartesian diagram
$$\CD\dx@>\dot m>>\dfg\\@Vm'VV @V\pi VV\\X@>m>>\fg\endCD$$
Then $m^*(\ck)$ is naturally an object of the bounded derived category of
constructible sheaves on $\dx$ with a $W$-action inherited from $\ck$ hence 
there is a natural $W$-action on the hypercohomology
$$H^j_c(X,m^*(\ck))=H^j_c(X,m'_!\dcl^*)=H^j_c(\dx,\dcl^*).$$
(We will often denote various local systems
obtained from $\dcl,\dcl^*$ by some natural construction again by
$\dcl,\dcl^*$.) 

\subhead 1.10 \endsubhead
If, in addition, $X$ has a given action of a closed connected subgroup $G'$ of
$G\tim\bc^*$ and $m$ is compatible with the $G'$-actions, and if $\Ga$ is a
smooth irreducible variety with a free $G'$-action, we can form the cartesian
diagram
$$\CD {}_\Ga\dx@>{}_\Ga\dot m>>{}_\Ga\dfg\\@V{}_\Ga m'VV @V{}_\Ga\pi VV\\
        {}_\Ga X     @>{}_\Ga m>>  {}_\Ga\fg  \endCD$$
where $Y\mto{}_\Ga Y$ is the functor from algebraic varieties with $G'$-action
to algebraic varieties given by $Y\mto G'\bsl(\Ga\tim Y)$.

The local system $\bc\bxt\dcl^*$ on $\Ga\tim\dfg$ is $G'$-equivariant hence it
descends canonically to a local system ${}_\Ga\dcl^*$ on ${}_\Ga\dfg$. Also,
$\bc\bxt\ck$ is (up to shift) a $G'$-equivariant perverse sheaf with $W$-action
on $\Ga\tim\fg$ hence it descends to a perverse sheaf (up to shift) ${}_\Ga\ck$
with $W$-action on ${}_\Ga\fg$. We have canonically 
${}_\Ga\ck=({}_\Ga\pi)_!({}_\Ga\dcl^*)$. Then $({}_\Ga m)^*({}_\Ga\ck)$ is 
naturally an object of the bounded derived category of constructible sheaves on
${}_\Ga\dx$, with a $W$-action inherited from ${}_\Ga\ck$ hence there is a 
natural $W$-action on the hypercohomology
$$H^{2d-j}_c({}_\Ga X,({}_\Ga m)^*({}_\Ga\ck))
=H^{2d-j}_c({}_\Ga X,({}_\Ga m')_!
(({}_\Ga\dot m)^*({}_\Ga\dcl^*)))=H^{2d-j}_c({}_\Ga\dx,{}_\Ga\dcl^*)$$
where $d=\dim\dx$. (We write ${}_\Ga\dcl^*$ instead of 
$({}_\Ga\dot m)^*({}_\Ga\dcl^*)$.) We can choose $\Ga$ so that $H^n(\Ga,\bc)=0$
for $n\in[1,m]$ where $m$ is large compared with $j$. Taking duals, we see that
$W$ acts naturally on the equivariant homology 
$$H_j^{G'}(\dx,\dcl)=H^{2d-j}_c({}_\Ga\dx,{}_\Ga\dcl^*)^*$$
\cite{\LI, 1.1}. 
This action is 
independent of the choice of $\Ga$.

\subhead 1.11\endsubhead
Let $\bs=S(\fh^*\opl\bc)=S(\fh^*)\ot\bc[\bor]$ where $S()$ denotes the 
symmetric algebra of a $\bc$-vector space and $\bor=(0,1)\in\fh^*\opl\bc$.

For any algebraic group $G'$ we write $H^*_{G'}$ instead of 
$H^*_{G'}(\text{\rm point},\bc)$ (equivariant cohomology). For any surjective 
homomorphism $G'@>>>G''$ of connected algebraic groups we have a canonical 
algebra homomorphism $H^*_{G'}@>>>H^*_{G''}$. (Using the identification 
\cite{\LI, 1.11(a)}, this is obtained by associating to a polynomial function 
$\ug''@>>>\bc$ its composition with the obvious map $\ug'@>>>\ug''$.) In 
particular, if $P\in\cp$, we have a canonical algebra 
homomorphism 
$$H^*_{P/P_{der}\tim\bc^*}=H^*_{\bap/\bap_{der}\tim\bc^*}@>>>
H^*_{\bap\tim\bc^*}.$$
Composing this with the algebra homomorphism
$H^*_{\bap\tim\bc^*}@>>>H^*_{\bap\tim\bc^*}(\bocp,\bc)$ (as in \cite{\LI, 1.7})
we obtain an algebra homomorphism 
$H^*_{P/P_{der}\tim\bc^*}@>>>H^*_{\bap\tim\bc^*}(\bocp,\bc)$. By 
\cite{\LI, 1.6, 1.4(e), 1.4(h)} we have canonically
$$H^*_{G\tim\bc^*}(\dfg,\bc)=H^*_{P\tim\bc^*}(\upi_P\i(\bocp+\uzp),\bc)=
H^*_{P\tim\bc^*}(\bocp,\bc)=H^*_{\bap\tim\bc^*}(\bocp,\bc).$$
We obtain an algebra homomorphism
$H^*_{P/P_{der}\tim\bc^*}@>>>H^*_{G\tim\bc^*}(\dfg,\bc)$. Using the canonical 
isomorphism $P/P_{der}@>\sim>>\bt$ we obtain an algebra homomorphism \linebreak
$H^*_{\bt\tim\bc^*}@>>>H^*_{G\tim\bc^*}(\dfg,\bc)$. This is in fact an algebra 
isomorphism (a reformulation of \cite{\LI, 4.2}). By \cite{\LI, 1.10} we have 
canonically $H^*_{\bt\tim\bc^*}=S(\fh^*\opl\bc)=\bs$. Thus we have an algebra 
isomorphism
$$\bs@>\sim>>H^*_{G\tim\bc^*}(\dfg,\bc).\tag a$$

Assume that $\tx$ is an algebraic variety with a given action of a closed 
connected subgroup $G'$ of $G\tim\bc^*$ and with a given morphism
$\tm:\tx@>>>\dfg$ compatible with the $G'$-actions. We write $\dcl$ instead of
$\tm^*\dcl$ (a local system on $\tx$). Then $H^{G'}_*(\tx,\dcl)$ is an 
$\bs$-module as follows. Let $\xi\in\bs$, let 
$\xi'\in H^*_{G\tim\bc^*}(\dfg,\bc)$ be the element that corresponds to $\xi$ 
under (a) and let $\xi''\in H^*_{G'}(\dfg,\bc)$ be the image of $\xi'$ under 
the homomorphism $H^*_{G\tim\bc^*}(\dfg,\bc)@>>>H^*_{G'}(\dfg,\bc)$ as in 
\cite{\LI, 1.4(f)}. We have $\tm^*(\xi'')\in H^*_{G'}(\tx,\bc)$. If
$z\in H^{G'}_*(\tx,\dcl)$ then $\xi z$ is defined as the product
$\tm^*(\xi'')\cdot z\in H_*^{G'}(\tx,\dcl)$ as in \cite{\LI, 1.7}.

\subhead 1.12\endsubhead 
Let $y\in\fg_N$. Let 
$$\cp_y=\{P\in\cp;y\in\up;\upi_P(y)\in\bocp+\uzp\}=
\{P\in\cp;y\in\up;\upi_P(y)\in\bocp\}.$$
(The second equality follows from the fact that $y$ is nilpotent.) The second 
projection identifies $\{y\}\do$ with $\cp_y$. Note that $\{y\}$ and $\{y\}\do$
are stable under $M(y)$ hence under $M^0(y)$. By 1.10 applied to $X=\{y\}$ and 
by 1.11 applied to $\tx=\{y\}\do$ we see that, if $G'$ is a closed connected 
subgroup of $M^0(y)$, then $H_*^{G'}(\cp_y,\dcl)$ has a natural $W$-action and
a natural $\bs$-action. It also has a natural $H^*_{G'}$-module structure 
\cite{\LI, 1.7}. Now there is a unique $\bh$-module structure on 
$H_*^{G'}(\cp_y,\dcl)$ such that $\bor\in\bh$ acts as $\bor\in\bs$, 
$\uxi\in\bh$ acts as $\xi\in\bs$ (for $\xi\in\fh^*$) and $s_i\in\bh$ acts as 
$s_i\in W$ (for $i\in I$). (In the special case when $G'=M^0(y)$, this follows
from \cite{\LI, 8.13}. The case when $G'$ is not necessarily $M^0(y)$ can be 
reduced to the special case using the isomorphism
$$H^*_{G'}\ot_{H^*_{M^0(y)}}H_*^{M^0(y)}(\cp_y,\dcl)@>\sim>>
H_*^{G'}(\cp_y,\dcl)\tag a$$
as in \cite{\LI, 7.5}; that result is applicable in view of \cite{\LI, 8.6}.)
The $\bh$-module structure commutes with the $H^*_{G'}$-module structure on 
$H_*^{G'}(\cp_y,\dcl)$. 

Now the finite group $\bam(y)=M(y)/M^0(y)$ acts on $H_*^{M^0(y)}(\cp_y,\dcl)$
by \cite{\LI, 1.9(a)}. This action commutes with the $\bh$-module structure and
is compatible with the $H^*_{M^0(y)}$-module structure where we regard
$H^*_{M^0(y)}$ as being endowed with the action of $\bam(y)$ given again by
\cite{\LI, 1.9(a)}.

\subhead 1.13\endsubhead 
Let $G'$ be a closed connected subgroup of $G\tim\bc^*$. Then 
$\un G'\sub\fg\opl\bc$. By \cite{\LI, 1.11(a)}, we may identify $H^*_{G'}$ with
the space of polynomials $f:\un G'@>>>\bc$ that are constant on the cosets by 
the nil-radical of $\un G'$ and are constant on the $\Ad$-orbits of $G'$. Let
$(\si,r)\in\un G'$ be a semisimple element.  Let $\cj^{G'}_{\si,r}$ be the 
maximal ideal of $H^*_{G'}$ consisting of all $f$ such that $f(\si,r)=0$. Let
$\bc_{\si,r}=H^*_{G'}/\cj^{G'}_{\si,r}$. (A one dimensional $\bc$-vector 
space.) 

Now assume that $G'\sub M^0(y)$. Then 
$\un G'\sub\{(x,r)\in\fg\opl\bc;[x,y]=2ry\}=\un{M^0(y)}$. In particular we have
$[\si,y]=2ry$.

Let $f_1:\un G'@>>>\bc$ be defined by $f_1(x,r)=r$. In the $H^*_{G'}$-module 
structure on $H_*^{G'}(\cp_y,\dcl)$, $f_1$ acts as multiplication by 
$\bor\in\bh$. We form
$$E_{y,\si,r}=\bc_{\si,r}\ot_{H^*_{G'}}H_*^{G'}(\cp_y,\dcl)=
H_*^{G'}(\cp_y,\dcl)/\cj^{G'}_{\si,r}H_*^{G'}(\cp_y,\dcl).$$
Then $E_{y,\si,r}$ inherits from $H_*^{G'}(\cp_y,\dcl)$ an $\bh$-module
structure in which $\bor\in\bh$ acts as multiplication by $r$ (since 
$f_1-r\in\cj^{G'}_{\si,r}$). By 1.12(a), $E_{y,\si,r}$ defined in terms of $G'$
is the same as $E_{y,\si,r}$ defined in terms of $M^0(y)$. For this reason we 
do not include $G'$ in the notation for $E_{y,\si,r}$.

Let $M(y,\si)=M(y)\cap(Z(\si)\tim\bc^*)$. Let $\bam(y,\si)$ be the group of 
connected components of $M(y,\si)$. The obvious map $\bam(y,\si)@>>>\bam(y)$ is
injective, since $M^0(y)\cap(Z(\si)\tim\bc^*)$ is connected. Clearly, the
restriction of the $\bam(y)$ action on $H^*_{M^0(y)}$ to $\bam(y,\si)$ leaves
$\cj^{M^0(y)}_{\si,r}$ stable hence the action of $\bam(y)$ on 
$H_*^{M^0(y)}(\cp_y,\dcl)$ induces an action of $\bam(y,\si)$ on 
$E_{y,\si,r}$. This action commutes with the $\bh$-module structure.

\subhead 1.14\endsubhead 
Let $\Irr\bam(y,\si)$ be a set of representatives for the isomorphism classes
of irreducible representations of $\bam(y,\si)$. For $\rh\in\Irr\bam(y,\si)$ 
let 
$$E_{y,\si,r,\rh}=\Hom_{\bam(y,\si)}(\rh,E_{y,\si,r}).$$ 
Let $\Irr_0\bam(y,\si)$ be the set of those $\rh\in\Irr\bam(y,\si)$ such that
$E_{y,\si,r,\rho}\ne 0$ or, equivalently (see \cite{\LI, 8.10}) such that $\rh$
appears in the restriction of the $\bam(y)$-module $H_*(\cp_y,\dcl)$ to
$\bam(y,\si)$. (Equivariant homology or cohomology in which the group is not
specified is understood to be with respect to the group $\{1\}$.)

\proclaim{Theorem 1.15}(a) Let $y,\si,r$ be as above; assume that $r\ne 0$. Let
$\rh\in\Irr_0\bam(y,\si)$. Then the $\bh$-module $E_{y,\si,r,\rh}$ has a unique
maximal submodule. Let $\bae_{y,\si,r,\rh}$ be the simple quotient of
$E_{y,\si,r,\rh}$.

(b) Let $r\in\bc^*$. The map $(y,\si,\rh)\mto\bae_{y,\si,r,\rh}$ establishes a
bijection between the set of all triples $(y,\si,\rh)$ with $y\in\fg_N$,
$\si\in\fg_{ss}$ with $[\si,y]=2ry$ and $\rh\in\Irr_0\bam(y,\si)$ (modulo
the natural action of $G$) and the set of isomorphism classes of simple 
$\bh$-modules in which $\bor$ acts as multiplication by $r$.
\endproclaim
The proof is given in 3.39, 3.41, 3.42. (A bijection as in (b) has already
been obtained in \cite{\LII} by other means since (a) was not known  in
\cite{\LII}.)

\subhead 1.16\endsubhead
Let $J\sub I$ and let $Q\in\cp_J$. Let $Q^1$ be a Levi subgroup of $Q$. Let 
$y\in\uq^1$ be nilpotent. Now $Q^1$ carries a cuspidal datum
$(\cp',\boc',\cl)$ analogous to that of $G$ (see 1.3). Here we identify
$Q^1=\baq$ via $\pi_Q$. Replacing $G$ by $Q^1$ in the definition of
$W,\fh,\bs,\bh,\cp_y,\dcl$ we get $W_J,\fh,\bs,\bh',\cp'_y,\dcl$. We use 
$P'\mto\pi_Q\i(P')$ to identify $\cp'$ with $\cp^*=\{P\in\cp;P\sub Q\}$ and 
$\cp'_y$ with $\cp^*_y=\{P\in\cp_y;P\sub Q\}$.

Let $C$ be a maximal torus of $M^0_{Q^1}(y)\sub M^0(y)$. Then 
$H^C_*(\cp^*_y,\dcl)$ is an $\bh'$-module (by 1.12 for $Q^1=\baq$ instead of 
$G$). Using the obvious algebra homomorphism $\bh'@>>>\bh$ (taking the 
generators of $\bh'$ to the corresponding generators of $\bh$) we can form the
$\bh$-module $\bh\ot_{\bh'}H^C_*(\cp^*_y,\dcl)$. Now the closed imbedding 
$j:\cp^*_y@>>>\cp_y$ induces a map 
$j_!:H^C_*(\cp^*_y,\dcl)@>>>H^C_*(\cp_y,\dcl)$, (see \cite{\LI, 1.4(b)}.) From
the definitions we see that $j_!$ is $\bh'$-linear hence it induces an
$\bh$-linear map
$$\bh\ot_{\bh'}H^C_*(\cp^*_y,\dcl)@>>>H^C_*(\cp_y,\dcl).\tag a$$
Let ${}_y\uuq=\cok(\ad(y):\uuq@>>>\uuq)$. Define $\ep:\un{M^0_{Q^1}(y)}@>>>\bc$
(recall that $\un{M^0_{Q^1}(y)}\sub\fg\opl\bc$) by
$$\ep(x,\la)=\det(\ad(x)-2\la:{}_y\uuq@>>>{}_y\uuq).$$
(For $(x,r)\in\un{M^0_{Q^1}(y)}$ we have $[x,y]=2ry$ hence 
$[\ad(x),\ad(y)]=2r\ad(y)$ hence $\ad(x):\uuq@>>>\uuq$ induces a map 
${}_y\uuq@>>>{}_y\uuq$ denoted again by $\ad(x)$.) The restriction of $\ep$ to
$\un C$ is denoted again by $\ep$. By the identification \cite{\LI, 1.11(a)} 
we may regard $\ep$ as an element of $H^*_C$. Applying 
$H^*_C[\ep\i]\ot_{H^*_C}$ to (a) we obtain an $H^*_C[\ep\i]$-linear map
$$H^*_C[\ep\i]\ot_{H^*_C}(\bh\ot_{\bh'}H^C_*(\cp^*_y,\dcl))@>>>
H^*_C[\ep\i]\ot_{H^*_C}H^C_*(\cp_y,\dcl).\tag b$$

\proclaim{Theorem 1.17 (Induction theorem)}The map 1.16(b) is an ($\bh$-linear)
isomorphism.
\endproclaim
The proof is given in Section 2 as an application of the "strong induction 
theorem" 2.16 which states the existence of an isomorphism similar to 1.16(b) 
but in which no elements of $H^*_C$ need to be inverted. 

\proclaim{Corollary 1.18} Assume that $(\si,r)$ is a semisimple element of
$\un{M^0_{Q^1}(y)}$ such that $\ep(\si,r)\ne 0$. Define $E'_{y,\si,r}$ like
$E_{y,\si,r}$ but in terms of $Q^1$ instead of $G$. Choose $C$ as in 1.16 such
that $(\si,r)\in\un C$. Then the map 1.16(b) induces an isomorphism of 
$\bh$-modules $\bh\ot_{\bh'}E'_{y,\si,r}@>\sim>>E_{y,\si,r}$.
\endproclaim

\subhead 1.19\endsubhead 
Let $\ci$ be the collection of all simple finite dimensional $\fg$-modules $V$
such that for any $P\in\cp$ there exists a $\up$-stable line $D_P$ in $V$ 
(necessarily unique). An equivalent condition is that, for any $P\in\cp$, 
$\{v\in V;\uup v=0\}$ is a line in $V$. (Clearly, $\{v\in V;\uup v=0\}$ is
$\up$-stable, hence if the second condition holds then the first condition 
holds. Conversely, let $P\in\cp$ and let $\fb$ be a Borel subalgebra of $\up$.
If $\{v\in V;\uup v=0\}$ is non-zero, then it is a $\up/\uup$-module with a 
unique line stable under $\fb/\uup$, a Borel subalgebra of $\up/\uup$, since 
such a line must be $\fb$-stable and $V$ is simple. It follows that 
$\{v\in V;\uup v=0\}$ is simple as a $\up/\uup$-module. Hence, if there exists
a $\up$-stable line in $V$ (necessarily contained in $\{v\in V;\uup v=0\}$) 
then that line must be equal to $\{v\in V;\uup v=0\}$ so that 
$\{v\in V;\uup v=0\}$ is a line.)

Let $V\in\ci$. Then $V$ defines an element $\xi_V\in\fh^*$ as follows. Let 
$P\in\cp$. Then $xv=u(x)v$ for $x\in\up,v\in D_P$ where $u:\up@>>>\bc$ is a Lie
algebra homomorphism. Then $u$ factors through a linear form on 
$\up/[\up,\up]=\fh$ denoted by $\xi_V$. It is independent of the choice of $P$.

\subhead 1.20\endsubhead 
{\it In the remainder of this section we assume that $G$ is semisimple.}

Let $r\in\bc^*$ and let $\ta:\bc@>>>\br$ be a homomorphism of abstract groups
such that $\ta(r)\ne 0$. Let $E$ be an $\bh$-module of finite dimension over 
$\bc$. We say that $E$ is {\it $\ta$-tempered} if for any $V\in\ci$, any 
eigenvalue $\la$ of $\un{\xi_V}$ on $E$ satisfies $\ta(\la)/\ta(r)\ge 0$. We 
say that $E$ is {\it $\ta$-square integrable} if for any $V\in\ci$, other than
$\bc$, any eigenvalue $\la$ of $\un{\xi_V}$ on $E$ satisfies 
$\ta(\la)/\ta(r)>0$.

The standard basis of the Lie algebra $\fs\fl_2(\bc)$ is denoted as follows:
$$e_0=\sm 0&1\\0&0\esm, h_0=\sm 1&0\\0&-1\esm, f_0=\sm 0&0\\1&0\esm.$$

\proclaim{Theorem 1.21} Assume that $r\ne 0$. Let $y,\si,\rh$ be as in 1.15(b).
The following three conditions are equivalent:

(i) $E_{y,\si,r,\rh}$ is $\ta$-tempered.

(ii) $\bae_{y,\si,r,\rh}$ is $\ta$-tempered.

(iii) There exists a homomorphism of Lie algebras $\ph:\fs\fl_2(\bc)@>>>\fg$ 
such that $y=\ph(e_0),[\si,\ph(h_0)]=0,[\si,\ph(f_0)]=-2r\ph(f_0)$ and any 
eigenvalue $\la$ of $\ad(\si-r\ph(h_0)):\fg@>>>\fg$ satisfies $\ta(\la)=0$.
\nl
If these conditions are satisfied, then $E_{y,\si,r,\rh}=\bae_{y,\si,r,\rh}$.
\endproclaim
The proof is given in 3.43.

\proclaim{Theorem 1.22} Assume that $r\ne 0$. Let $y,\si,\rh$ be as in 1.15(b).
The following five conditions are equivalent:

(i) $y,\si$ are not contained in a Levi subalgebra of a proper parabolic 
subalgebra of $\fg$.

(ii) There exists a homomorphism of Lie algebras $\ph:\fs\fl_2(\bc)@>>>\fg$ 
such that $y=\ph(e_0),\si=r\ph(h_0)$; moreover, $y$ is distinguished.

(iii) $\bae_{y,\si,r,\rh}$ is $\ta$-square integrable.

(iv) $E_{y,\si,r,\rh}$ is $\ta$-square integrable.

(v) For any $V\in\ci$, other than $\bc$, any eigenvalue of $r\i\un{\xi_V}$ on
$E_{y,\si,r,\rh}$ is an integer $\ge 1$.
\nl
If these conditions are satisfied, then $E_{y,\si,r,\rh}=\bae_{y,\si,r,\rh}$.
\endproclaim
The proof is given in 3.44.

\head 2. A strong induction theorem and a proof of Theorem 1.17\endhead
\subhead 2.1\endsubhead 
In this section we place ourselves in the setup of 1.16. Thus, 
$$J,Q,Q^1,\bh',y,\cp',\cp'_y,\cp^*,\cp^*_y,C$$ 
are defined. We set $\fq=\uq,\fq^1=\uq^1,\fn=\uuq$. We can find a Lie algebra 
homomorphism $\ph:\fs\fl_2(\bc)@>>>\fq^1$ such that $\ph(e_0)=y$ and such that
$C$ is a maximal torus of 
$$\{(g,\la)\in Q^1\tim\bc^*;\Ad(g)\ph(e_0)=\la^2\ph(e_0),
\Ad(g)\ph(f_0)=\la^{-2}\ph(f_0)\}.$$
For any $P\in\cp$ we set 
$$P^!=(P\cap Q)U_Q$$ 
(a parabolic subgroup of $Q$).

Let $f:\fq@>>>\fq^1$ be the projection of $\fq=\fq^1\opl\fn$ onto $\fq^1$. Let
$$W_*=\{w\in W; w \text{ has minimal length in } wW_J\}.$$
There are only finitely many orbits for the conjugation action of $Q$ on $\cp$.
A $Q$-orbit $\co$ on $\cp$ is said to be {\it good} if any $P\in\co$ has some
Levi subgroup that is contained in $Q$. For such $\co$ and for $P\in\co$ we 
have $P^!\in\cp^*$ and the $G$-orbit of $(P,P^!)$ in $\cp\tim\cp$ is good (see
1.5) and indexed by an element $w\in W_*$. Moreover, $\co\mto w$ is a well 
defined bijection between the set of good $Q$-orbits on $\cp$ and $W_*$. We 
denote by $\boo(w)$ the $Q$-orbit on $\cp$ corresponding to $w\in W_*$. Note 
that $P\mto P^!$ is a $Q$-equivariant morphism $\boo(w)@>>>\cp^*$. 

If $X$ is a subvariety of $\fg$ then $\dx=\{(z,P)\in\dfg;z\in X\}$ (see 1.9) is
a subvariety of $\dfg$. For any subvariety $F$ on $\cp$ we set
$$\dx_F=\{(z,P)\in\dx;P\in F\}.$$ 
For $w\in W_*$ we will often write $\dx_w$ instead of $\dx_{\boo(w)}$. Define 
$$\tf:\dq_w@>>>(q^1)\dot{}_1=\dq^1_1, (z,P)\mto(f(z),P^!).$$

Let $S=Z_{Q^1}^0$. Let $w\in W_*$. Let $\boo(w)^S$ be the fixed point set of 
the conjugation action of $S$ on $\boo(w)$. The properties (a),(b) below are 
easily verified.

(a) The map $\boo(w)^S@>>>\cp^*$ given by $P@>>>P^!$ is an isomorphism. 

(b) Let $z\in\fq$ and let $P\in\boo(w)^S$. Then we have $(z,P)\in\dfg$ if and
only if $(z,P^!)\in\dfg$.
\nl
The fixed point set of the $S$-action on $\dq_w$ (conjugation on both factors)
is $\dq^1_w$. The restriction of $\tf$ defines an isomorphism 
$\dq^1_w@>\sim>>\dq^1_1$. 

We choose a homomorphism of algebraic groups $\chi:\bc^*@>>>S$ such that
$\la\mto\Ad(\chi(\la))$ has weights $>0$ on $\fn$. We define a $\bc^*$-action
on $\dq_w$ by 
$$\la:(z,P)\mto(\Ad(\chi(\la))z,\chi(\la)P\chi(\la)\i).$$
Then $\tf:\dq_w@>>>\dq^1_1$ is $\bc^*$-equivariant where $\bc^*$ acts on 
$\dq^1_1$ trivially. Let $n=\dim\bap-\dim\bocp-\dim Z_{\bap}$ for any 
$P\in\cp$.

\proclaim{Lemma 2.2} (a) $\dq_w$ is a smooth variety of pure dimension
$\dim Q-n$.

(b) Let $(z',P')\in\dq_w$; define $(z,P)\in\dq^1_w$ by $\tf(z',P')=\tf(z,P)$. 
Then \linebreak
$\lim_{\la\to 0}\la(z',P')$ exists in $\dq_w$ and equals $(z,P)$.

(c) The fixed point set of the $\bc^*$-action on $\dq_w$ coincides with the
fixed point set of the $S$-action on $\dq_w$.
\endproclaim
Consider the fibration $pr_2:\dq_w@>>>\boo(w)$. Let $P\in\boo(w)$. Let $P^1$ be
a Levi subgroup of $P$ that is contained in $Q$. Let $c^1$ be the nilpotent 
orbit in $\un P^1$ corresponding to $\bocp$ under $P^1@>\sim>>\bap$. Then 
$pr_2\i(P)$ may be identified with 

$(c^1+Z_{P^1}+\uup)\cap\fq=c^1+Z_{P^1}+(\uup\cap\fq)$
\nl
and this is smooth irreducible of dimension 

$-n+\dim P^1+\dim(U_P\cap Q)=-n+\dim(P\cap Q)$. 
\nl
Now $\boo(w)$ is smooth,
irreducible of dimension $\dim Q-\dim(P\cap Q)$ and (a) follows.

We prove (b). We have $z'=z+x$ with $x\in\fn$. We have $P'{}^!=P^!$ and $P',P$ 
are in $\boo(w)$ hence $P'=uPu\i$ for some $u\in U_Q$. Thus,
$$\align&\la(z',P')=(\Ad(\chi(\la))(z+x),\chi(\la)uPu\i\chi(\la)\i)\\&
=(z+\Ad(\chi(\la))x,\chi(\la)u\chi(\la\i)P\chi(\la)u\i\chi(\la)\i).\endalign$$
Now $\lim_{\la\to 0}\chi(\la)u\chi(\la\i)=1$ and
$\lim_{\la\to 0}\Ad(\chi(\la))x=0$ by the choice of $\chi$. Hence 
$\lim_{\la\to 0}\la(z',P')=(z,P)$. This proves (b). Now (c) follows immediately
from (b).

\subhead 2.3\endsubhead
Let $A=y+\fz_\fq(\ph(f_0)),\quad A^1=A\cap\fq^1.$

\proclaim{Lemma 2.4} Let $w\in W_*$. 

(a) $\da_w$ is smooth of pure dimension $-n+\dim\fz_\fq(\ph(f_0))$. 

(b) $\da^1_1$ is smooth of pure dimension $-n+\dim\fz_{\fq^1}(\ph(f_0))$. 

(c) The map $\da_w@>>>\da^1_1$, $(z,P)\mto(f(z),P^!)$ is an affine space bundle
with fibres of dimension $\dim\fz_\fn(\ph(f_0))$.
\endproclaim
By \cite{\LB, 6.9},

(d) the map $Q\tim A@>>>\fq$ (given by the adjoint action of $Q$) is smooth 
with fibres of pure dimension $\dim\fz_\fq(\ph(f_0))$.
\nl
We have a cartesian diagram
$$\CD \da_w@>pr_1>>A\\ @VVV @VVV\\ \dq_w @>pr_1>>\fq\endCD$$
where the vertical maps are the inclusions. This induces a cartesian diagram
$$\CD Q\tim\da_w@>>>Q\tim A\\ @VVV @VVV \\ \dq_w @>>>\fq\endCD$$
where the vertical maps are given by the obvious action of $Q$. Using (d), it
follows that $Q\tim\da_w@>>>\dq_w$ (adjoint action) is smooth with fibres of
pure dimension $\dim\fz_\fq(\ph(f_0))$. Since $\dq_w$ is smooth of pure
dimension $\dim Q-n$ (see 2.2(a)), it follows that $Q\tim\da_w$ is smooth of
pure dimension $\dim Q-n+\dim\fz_\fq(\ph(f_0))$. Hence $\da_w$ is as in (a).
Now $\da^1_1$ is the same as $\da_1$ (where $\fg,\fq$ are replaced by 
$\fq^1,\fq^1$). Hence (b) follows from (a).

We prove (c). Note that $\da_w$ is a closed subset of $\dq_w$, stable under the
$S$-action (hence under the $\bc^*$-action) on $\dq_w$ (as in 2.1). By 2.2(c),
the fixed point set of the $\bc^*$-action on $\da_w$ is the same as the fixed 
point set of the $S$-action of $\da_w$ that is,

$\da_w^{\bc^*}=\{(z,P)\in\dfg;z\in A^1,P\in\boo(w)^S\}$
\nl
and the map $\da_w@>>>\da^1_1$ (restriction of $\tf:\dq_w@>>>\dq^1_1$) 
restricts to an isomorphism $\da_w^{\bc^*}@>\sim>>\da^1_1$. Using 2.2(b) and 
the fact that $\da_w$ is closed in $\dq_w$, we see that

(e) for any $(z',P')\in\da_w$, $\lim_{\la\to 0}\la(z',P')$ exists in $\da_w$ 
and belongs to $\da_w^{\bc^*}$.
\nl
Let $Z$ be a connected component of $\da_w$. Then
$Z^{\bc^*}=Z\cap\da_w^{\bc^*}$ may be identified with $\tf(Z)$ hence is 
connected (a connected component of $\da_w^{\bc^*}\cong\da^1_1$). Let $Z'$ be a
smooth $\bc^*$-equivariant projective compactification of $Z$. We have 
$Z^{\bc^*}\sub Z'{}^{\bc^*}$; let $Z'_1$ be the connected component of 
$Z'{}^{\bc^*}$ that contains $Z^{\bc^*}$. Let $Z'_2$ be the set of all 
$x\in Z'$ such that $\lim_{\la\to 0}\la x$ exists in $Z'$ and belongs to 
$Z'_1$. From (e) we see that $Z\sub Z'_2$; hence $Z'_2$ is dense in $Z'$. By a
known result of Bialynicki-Birula, $Z'_2$ is locally closed in $Z'$ (hence
open) and the map $\tf':Z'_2@>>>Z'_1$ given by $x\mto\lim_{\la\to 0}\la x$ is
an affine space bundle. Let $x\in Z^{\bc^*}$, $Z^x=\{x'\in Z;\tf(x')=x\}$. We
have

$\dim\tf'{}\i(x)=\dim Z'_2-\dim Z'_1=\dim Z-\dim Z'_1\le\dim Z-\dim Z^{\bc^*}
\le\dim Z^x$.
\nl
Since $Z^x\sub\tf'{}\i(x)$ and $\tf'{}\i(x)$ is irreducible, we see that we
must have $Z^x=\tf'{}\i(x)$ and $\dim Z^{\bc^*}=\dim Z'_1$. Thus we have a 
cartesian diagram
$$\CD Z@>>>Z'_2\\@VVV @VVV\\Z^{\bc^*}@>>>Z'_1\endCD$$
(the horizontal maps are inclusions and the vertical maps are $\tf,\tf'$). It
follows that $\tf:Z@>>>Z^{\bc^*}$ is an affine space bundle with fibres of
dimension 
$$\align&
\dim Z'_2-\dim Z'_1=\dim Z-\dim Z^{\bc^*}=\dim\da_w-\dim\da_w^{\bc^*}=\\&
-n+\dim\fz_\fq(\ph(f_0))-(-n+\dim\fz_{\fq^1}(\ph(f_0)))=\dim\fz_\fn(\ph(f_0)).
\endalign$$
(c) follows. The lemma is proved.

\subhead 2.5\endsubhead
Let 
$$\ca=y+\fz_\fn(\ph(f_0)),\quad\ca^1=\ca\cap\fq^1=\{y\},$$
$$\ca'=y+\fz_\fn(\ph(f_0))+\ft,\quad\ca'{}^1=\ca'\cap\fq^1=y+\ft,$$
$$\ca''=y+\fz_\fn(\ph(f_0))+\ft_r,\quad \ca''{}^1=\ca''\cap\fq^1=y+\ft_r,$$
where $\ft=\fz_{\fq^1}$ and $\ft_r=\{x\in\ft;\fzg(x)=\fq^1\}$. We have 
$\ca\sub\ca'\supset\ca'',\ca^1\sub\ca'{}^1\supset\ca''{}^1$. Hence
$\dca\sub\dca'\supset\dca'',\dca^1_1\sub\dca'{}^1_1\supset\dca''{}^1_1$. We 
have 
$$\dca^1_1=\{y\}\tim\cp^*_y,\quad\dca'{}^1_1=(y+\ft)\tim\cp^*_y,\quad
\dca''{}^1_1=(y+\ft_r)\tim\cp^*_y.\tag a$$
From the definitions we see that, for $w\in W_*$, we have cartesian diagrams
$$\CD\dca_w@>>>\da\\@VVV @VVV\\ \dca^1_1@>>>\da^1_1\endCD$$
$$\CD\dca'_w@>>>\da\\@VVV @VVV\\ \dca'{}^1_1@>>>\da^1_1\endCD$$
$$\CD\dca''_w@>>>\da\\@VVV @VVV\\ \dca''{}^1_1@>>>\da^1_1\endCD$$
where the horizontal maps are the obvious inclusions and the vertical maps are
defined by $(z,P)\mto(f(z),P^!)$. Using this and 2.4(c) we see that 

(b) {\it the maps
$$\dca_w@>>>\dca^1_1,\quad\dca'_w@>>>\dca'{}^1_1,\quad\dca''_w@>>>\dca''{}^1_1
$$
defined by $(z,P)\mto(f(z),P^!)$ are affine space bundles with fibres of 
dimension} $\dim\fz_\fn(\ph(f_0))$.

\proclaim{Lemma 2.6} Let $E$ be a connected algebraic group, let $U$ be a
closed normal unipotent subgroup of $E$ and let $\ct$ be a torus in $E$. Let 
$e,h,e'$ be elements of $\un E$ such that $[h,e]=2e,[h,e']=-2e',[e,e']=h$ and 
such that $[d,e]=[d,e']=[d,h]=0$ for any $d\in\un\ct$. Let $\fu=\un U$. Define
$\Phi:U\tim\fz_\fu(e')\tim\un\ct@>>>\un\ct\tim\fu$ by 
$\Phi(u,x,d)=(d,\Ad(u)(e+x+d)-e-d)$. Then $\Phi$ is an affine space bundle with
fibres isomorphic to $\fz_\fu(e)$.
\endproclaim
Let $Z=Z_U$. We show that there exists a morphism $d\mto L'_d$, 
$\un\ct@>>>\Hom(\uz,\uz)$ and a morphism $d\mto L''_d$, 
$\un\ct@>>>\Hom(\uz,\ker(\ad(e'):\uz@>>>\uz))$ such that 
$$z=[L'_d(z),e+d]+L''_d(z)\tag a$$
for any $z\in\uz,d\in\un\ct$. We can find a direct sum decomposition 
$\uz=\opl_{k=1}^NV_k$ where $V_k$ are vector subspaces of $\uz$ and linear
forms $a_1,\dots,a_N$ on $\un\ct$ such that for any $k$, $V_k$ is a simple
$\fs\fl_2(\bc)$-submodule of $\uz$ under $\ad(e),\ad(h),\ad(e')$ and
$[d,v]=a_k(d)v$ for any $d\in\un\ct,v\in V_k$. 

Let $b_{0,k},b_{1,k},\dots,b_{n_k,k}$ be a basis of $V_k$ such that
$$\ad(e')b_{0,k}=0,\ad(e)b_{0,k}=-b_{1,k},\ad(e)b_{1,k}=-b_{2,k},\dots,
\ad(e)b_{n_k-1,k}=-b_{n_k,k}.$$
For $s<0$ we set $b_{s,k}=0$. For $s\in[0,n_k]$ we set
$$\align&L'_d(b_{s,k})=b_{s-1,k}+a_k(d)b_{s-2,k}+a_k(d)^2b_{s-3,k}+\dots
+a_k(d)^{s-1}b_{0,k},\\& L'_d(b_{s,k})=a_k(d)^sb_{0,k}.\endalign$$
This defines uniquely $L'_d,L''_d$. It is clear that (a) holds.

Next we construct for any $d\in\un\ct$ an isomorphism

(b) 
$p:\{(z',z)\in\uz\tim\fz_{\uz}(e');[z',e+d]+z=0\}@>\sim>>\fz_{\uz}(e)$.
\nl
This is by definition a direct sum over $k$ of isomorphisms
$$\align&p_k:\{(v',v)\in V_k\tim\ker(\ad(e'):V_k@>>>V_k);-\ad(e+d)(v')+v=0\}
\\&@>\sim>>\ker(\ad(e'):V_k@>>>V_k)\endalign$$
given by
$p_k(c_0b_{0,k}+c_1b_{1,k}+\dots+c_{n_k}b_{n_k,k},c'b_0)=c_{n-k}b_{n-k,k}$.

We prove the lemma by induction on $\dim U$. If $\dim U=0$, the result is 
trivial. Hence we may assume that $\dim U>0$ and that the result is true when 
$E,U$ are replaced by $\bae=E/Z,\bau=U/Z$, $\ct$ is replaced by the image
$\bct$ of $\ct$ in $\bae$ and $e,h,e'$ are replaced by their images 
$\bar e,\bar h,\bar e'$ under $\un E@>>>\un{\bae}$. (We have $\dim Z>0$.) Let
$\bfu=\un{\bau}$. The obvious map $\fu@>>>\bfu$ may be regarded as a surjective
map of $\fs\fl_2(\bc)$-modules; by the complete reducibility of such modules, 
this map admits a cross section as an $\fs\fl_2(\bc)$-module; in particular, 
the induced map $\fz_\fu(e')@>>>\fz_\fu(\bar e')$ is surjective. Let 
$\bax\mto\tbx$ be a linear cross section for this last linear map. Let
$\bar u\mto\tbu$ be an algebraic cross section $\bau@>>>U$ for the obvious map
$U@>>>\bau$. 

Let $(d,\xi)\in\un\ct\tim\fu$; let $\bar\xi$ be the image of $\xi$ under 
$\fu@>>>\bfu$ and let $\bar d$ be the image of $d$ under the isomorphism
$\ct@>\sim>>\bct$ induced by $E@>>>\bae$. We have
$$\Phi\i(d,\xi)\cong\{(u,x)\in U\tim\fz_\fu(e');\Ad(u)(e+x+d)=e+d+\xi\}.$$
By the induction hypothesis,
$$X=:\{(\bar u,\bax)\in\bau\tim\fz_{\bfu}(\bar e');
\Ad(\bar u)(\bar e+\bax+\bar d)=\bar e+\bar d+\bar\xi\}$$
is an affine space isomorphic to $\fz_{\bfu}(\bar e)$. We have an obvious map
$\Psi:\Phi\i(\xi)@>>>X$. Its fibre at $(\bar u,\bax)\in X$ is the set of all
$(u,x)\in U\tim\fz_\fu(e')\tim\ct$ such that $\Ad(u)(e+x+d)=e+d+\xi$ and 
$u=\tbu\ze,x=\tbx+z$ for some $\ze\in Z,z\in\fz_{\uz}(f)$. Note that 
$\Ad(\tbu)(e+\tbx+d)=e+d+\xi+z_0$ for some $z_0\in\uz$. The equation 
$\Ad(u)(e+x+d)=e+d+\xi$ can be written as $\Ad(\tbu\ze)(e+\tbx+z+d)=e+d+\xi$, 
or as $\Ad(\ze)(e+d+\xi+z_0)+z=e+d+\xi$, or as $\Ad(\ze)(e+d)+z=e+d-z_0$.
Setting $\ze=\exp(z'),z'\in\uz$, we see that the fibre of $\Psi$ at 
$(\bar u,\bax)$ may be identified with

(c) $\{(z',z)\in\uz\tim\fz_{\uz}(e');\Ad(\exp(z'))(e+d)+z=e+d-z_0\}$.
\nl
Since $[\uz,e+d]\in\uz$, we have $[z',[z',e+d]]=0$ for $z'\in\uz$ and (c)
becomes 

$\{(z',z)\in\uz\tim\fz_{\uz}(e');[z',e+d]+z=-z_0\}$
\nl
or

$\{(z',z)\in\uz\tim\fz_{\uz}(e');[z',e+d]+z+[L'_d(z_0),e+d]+L''_d(z_0)=0\}
$.
\nl
By the substitution $\ti z'=z'+L'_d(z_0),\ti z=z+L''_d(z_0)$, this becomes

$\{(\ti z',\ti z)\in\uz\tim\fz_{\uz}(e');[\ti z',e+d+\ti z=0\}$.
\nl
By the isomorphism (b) this is identified with the vector space $\fz_{\uz}(e)$.
We see that $\Phi\i(\xi)$ is a vector bundle over $X$. Since $X$ is an affine 
space, this vector bundle must be trivial (Quillen-Suslin) and therefore 
$\Phi\i(\xi)$ is itself an affine space isomorphic to 
$\fz_{\uz}(e)\tim\fz_{\bfu}(\bar e)\cong\fz_\fu(e)$.
The lemma is proved.

\proclaim{Lemma 2.7} Let $\ft'$ be the Lie algebra of a torus contained in 
$Z_{Q^1}$. Let $X=y+\fz_\fn(\ph(f_0))+\ft'$. Let $\co$ be a $Q$-orbit on $\cp$.
Assume that $\co$ is not good. Then $H^j_c(\dx_\co,\dcl^*)=0$ for any 
$j\in\bz$.
\endproclaim
In this proof all local systems are deduced from $\cl^*$ and we omit them from
the notation.

The assignment $P\mto P^!$ defines a morphism $\pi:\co@>>>\cp'$ where $\cp'$ is
a conjugacy class of parabolic subgroups of $Q$. The fibres of $\pi$ are
exactly the $U_Q$-orbits on $\co$. It is enough to show that
$H^i_c(\dx_F,)=0$ for any fibre $F$ of $\pi$ and any $i\in\bz$, where
$\dx_F=\{(z,P)\in\dx; P\in F\}$ (see 2.1). We fix $P\in F$; let $\fp=\up$. Let
$\cy=\{(z,u)\in\fg\tim U_Q;(Ad(u)z,P)\in\dfg,z\in X\}$. Define $\cy@>>>\dx_F$ 
by $(z,u)\mto(z,u\i Pu)$. This is a fibration with fibres isomorphic to 
$U_Q\cap P$. It is then enough to show that $H^i_c(\cy,)=0$ for any $i\in\bz$. 
Setting $z-y=x+d$ where $x\in\fz_\fn(\ph(f_0)),d\in\ft'$ and 
$\ti\bocp=\upi_P\i(\bocp+\uzp)$, we identify $\cy$ with 

$\{(x,d,u)\in\fz_\fn(\ph(f_0))\tim\ft'\tim U_Q;\Ad(u)(y+x+d)\in\ti\bocp\}$.
\nl
This maps to

$\cy'=\{(d,\nu)\in\ft'\tim\fn;y+d+\nu\in\ti\bocp\}$
\nl
by $(x,d,u)\mto\nu=(d,\Ad(u)(y+x+d)-y-d)$; this is an affine space bundle by 
2.6 applied to $E=Q,U=U_Q,\un\ct=\ft',e=y,h=\ph(h_0),e'=\ph(f_0)$. Hence it is
enough to show that $H^i_c(Y',)=0$. For any $d\in\ft'$ let

$\cy'_d=\{\nu\in\fn;y+d+\nu\in\ti\bocp\}$
\nl
be the fibre at $d$ of $pr_1:\cy'@>>>\ft'$. By the Leray spectral sequence for
$pr_1$, it is enough to show that $H^i_c(\cy'_d,)=0$ for any $d\in\ft'$. If 
$y+d\notin\fn+\fp$ then $\cy'_d=\emp$ and there is nothing to prove. Thus we 
may assume that $y+d+\nu_0\in\fp$ for some $\nu_0\in\fn$. Setting 
$\xi=y+d+\nu_0\in\fp\cap\fq,\nu'=\nu-\nu_0$, we may identify $\cy'_d$ with
$$\ti{\cy}'=\{\nu'\in\fn;\xi+\nu'\in\ti\bocp\}=
\{\nu'\in\fn\cap\fp;\xi+\nu'\in\ti\bocp\}.$$
Let $R=(P\cap Q)U_P$ (a proper parabolic subgroup of $P$ since $\co$ is not 
good). Let $\bar R$ be the image of $R$ under $P\mto\bap$ (a proper parabolic
subgroup of $\bap$). The nil-radical of $\un R$ is $\fn\cap\fp+\uup$. Hence the
nil-radical of $\ubar$ is
$$\fn_1=(\fn\cap\fp+\uup)/\uup=(\fn\cap\fp)/(\fn\cap\uup).$$
Let $k:\fn\cap\fp@>>>\fn_1$ be the canonical map. Let $\xi_1$ be the image of 
$\xi$ under $\fp@>>>\ubap$. We have $\xi\in\fq$ hence $\xi_1\in\ubar$. Let
$$\cy''=\{\mu\in\fn_1;\xi_1+\mu\in\bocp+\uzp\}.$$
We have a cartesian diagram
$$\CD\ti{\cy}'@>>>\fn\\@Vk'VV @VkVV\\ \cy''@>>>\fn_1\endCD$$
where the horizontal maps are the inclusions and $k'$ is induced by $k$. Since
$k$ is an affine space bundle with fibres isomorphic to $\fn\cap\uup$, so is
$k'$. It is therefore enough to show that $H^i_c(\cy'',)=0$. This follows
directly from the fact that $\cl^*$ is a cuspidal local system since $\fn_1$ is
the nil-radical of the proper parabolic algebra $\ubar$ of $\ubap$ and
$\xi_1\in\ubar$. The lemma is proved.

\proclaim{Lemma 2.8} Let $\de=\dim\fz_\fn(\ph(f_0)), b=\dim\ft$.

(a) We have $H^j_c(\dca',\dcl^*)=0,H^j_c(\dca'_1,\dcl^*)=0$ for odd $j$. For
any $j$ we have $\dim H^j_c(\dca',\dcl^*)=
\sha(W_*)\dim H^{j-2\de-2b}_c(\cp^*_y,\dcl^*)$.

(b) We have $H^j_c(\dca,\dcl^*)=0,H^j_c(\dca_1,\dcl^*)=0$ for odd $j$. For
any $j$ we have $\dim H^j_c(\dca,\dcl^*)
=\sha(W_*)\dim H^{j-2\de}_c(\cp^*_y,\dcl^*)$.

(c) Let $D=\dim\dca,D'=\dim\dca'$. For any $j$ we have
$$\dim H_j^C(\dca',\dcl)=\dim H_{j+2D-2D'+2b}^C(\dca,\dcl).$$
\endproclaim
We can arrange the $Q$-orbits on $\cp$ in a sequence
$\co_1,\co_2,\dots,\co_n$ so that $R_m=\co_1\cup\co_2\cup\dots\cup\co_m$ is 
closed in $\cp$ for any $m\in[1,n]$. We set $R_0=\emp$. For any $m$ we have 
$$H^j_c(\dca'_{\co_m},\dcl^*)=0 \text{ for odd } j.\tag d$$
(Indeed, if $\co_m$ is not good, this follows from 2.7; if $\co_m$ is good 
then, using 2.5(a),(b) we see that it is enough to show that 
$H^j_c(\cp^*_y,\dcl^*)=0$ for odd $j$. This follows from \cite{\LI, 8.6}.)
Using induction on $m$ we deduce that
$$H^j_c(\dca'_{R_m},\dcl^*)=0 \text{ for odd } j.\tag e$$
for any $m$. Taking $m=n$ we obtain the first sentence of (a). For $m\in[1,n]$
we have a cohomology exact sequence 
$$0@>>>H^j_c(\dca'_{\co_m},\dcl^*)@>>>
H^j_c(\dca'_{R_m},\dcl^*)@>>>H^j_c(\dca'_{R_{m-1}},\dcl^*)@>>>0\tag f$$
(we use (c),(d)). Using induction on $m$ it follows that

$\dim H^j_c(\dca'_{R_m},\dcl^*)=\sum_{m'\in[1,m]}
\dim H^j_c(\dca'_{\co_{m'}},\dcl^*)$.
\nl
Taking $m=n$ we obtain

$\dim H^j_c(\dca',\dcl^*)=\sum_{m'\in[1,n]}\dim H^j_c(\dca'_{\co_{m'}},\dcl^*)
$.
\nl
Using 2.7 and 2.5(a),(b) we obtain the second sentence in (a).

The proof of (b) is entirely similar (it is again based on 2.7 and 2.5(a),(b)).

We prove (c). The last equality in (c), in the non-equivariant case (that is 
the case where $C$ is replaced by $\{1\}$) follows immediately from the 
equation $\dim H^j_c(\dca',\dcl^*)=\dim H^{j-2b}_c(\dca,\dcl^*)$ (see 
(a),(b)). The case where $C$ is present can be deduced from the non-equivariant
case using the existence of (non canonical) isomorphisms of graded vector 
spaces
$$H_*^C(\dca',\dcl)\cong H^*_C\ot H_*(\dca',\dcl),\quad
H_*^C(\dca,\dcl)\cong H^*_C\ot H_*(\dca,\dcl)$$
which follows from \cite{\LI, 7.2(a)} (which is applicable, by (a),(b)). The 
lemma is proved.

(g) {\it Remark.} The following six conditions are equivalent:

(1) $H^*_c(\cp_y,\dcl^*)\ne 0$;

(2) $H^*_c(\cp^*_y,\dcl^*)\ne 0$;

(3) $H^*_c(\dca,\dcl^*)\ne 0$;

(4) $H^*_c(\dca_1,\dcl^*)\ne 0$;

(5) $H^*_c(\dca',\dcl^*)\ne 0$;

(6) $H^*_c(\dca'_1,\dcl^*)\ne 0$.

Indeed, we have $(5)\lra(2)$ by (a),$(3)\lra(2)$ by (b), $(6)\lra(2)$ and
$(4)\lra(2)$ by 2.5(a),(b). It remains to show that $(1)\lra(3)$. Using 
\cite{\LI, 7.2} (which is applicable by the odd vanishing (b) and 
\cite{\LI, 8.6}) we see that $(1)\lra(1')$ and $(3)\lra(3')$ where

(1') $H^C_*(\cp_y,\dcl)\ne 0$;

(3') $H^C_*(\dca,\dcl)\ne 0$.

It remains to show that $(1')\lra(3')$. Let $\cp_y^C,\dca^C,\fn^C$ be the fixed
point sets of the $C$-action on $\cp_y,\dca,\fn$. Since $C$ contains 
$Z_{Q^1}^0$ we have $\fn^C=\{0\}$ hence $\cp_y^C=\dca^C$. By the localization 
theorem \cite{\LII, 4.4} (which is applicable by \cite{\LI, 8.6} and (b)), 
the canonical $H^*_C$-linear maps $H^C_*(\cp_y^C,\dcl)@>>>H^C_*(\cp_y,\dcl)$,
$H^C_*(\dca^C,\dcl)@>>>H^C_*(\dca,\dcl)$ become isomorphisms after the scalars
are extended to the quotient field of $H^*_C$. Hence the canonical 
$H^*_C$-linear map $H^C_*(\cp_y,\dcl)@>>>H^C_*(\dca,\dcl)$ becomes an
isomorphism after the scalars are extended to the quotient field of $H^*_C$.
Since the $H^*_C$-modules $H^C_*(\cp_y,\dcl),H^C_*(\dca,\dcl)$ are finitely 
generated, projective (by \cite{\LI, 7.2} which is applicable by 
\cite{\LI, 8.6} and (b)) it follows that $(1')\lra(3')$. 

The previous argument shows also that $H^C_*(\cp_y,\dcl)@>>>H^C_*(\dca,\dcl)$ 
is injective.

\medpagebreak

For $x\in\fg$ let $x_s$ be the semisimple part of $x$.

\proclaim{Lemma 2.9}  Let $Y=\{x\in\fq;\fzg(x_s)\in\fq\}$.

(a) $Y$ is an open dense subset of $\fq$.

(b) Let $z$ be an element in the image of $pr_1:\dy@>>>Y$. There exists a Levi
subgroup $L$ of $Q$ such that the following holds: for any $P\in\cp$ such that
$(z,P)\in\dy$ we have $Z^0_L\sub P$. 

(c) We have $\dy=\cup_{w\in W_*}\dy_w$.

(d) For any $w\in W_*$, there is a well defined isomorphism of algebraic
varieties $f_w:\dy_w@>\sim>>\dy_1$ given by $(x,P)\mto(x,P^!)$. 

(e) For any $w\in W_*$, $\dy_w$ is open and closed in $\dy$.
\endproclaim
We prove (a). Let $\cs$ be the variety of semisimple classes of $\fq/\fn$. Then
$Y$ is the inverse image under $\fq@>>>\fq/\fn@>>>\cs$ (composition of 
canonical maps) of a subset $\bar Y$ of $\cs$. Let $\fc$ be a Cartan subalgebra
of $\fq^1$. It is enough to show that $\bar Y$ is non-empty, open in $\cs$ or 
that the inverse image $\ti Y$ of $\bar Y$ under the canonical open map 
$\fc@>>>\cs$ is non-empty open in $\fc$. Now
$$\ti Y=Y\cap\fc=\{x\in\fc;\fzg(x)\sub\fq\}=\{x\in c;\fzg(x)\sub\fq^1\}.$$
Let $R_0$ (resp. $R'_0$ be the set of roots of $\fg$ (resp. of $\fq^1$) with
respect to $\fc$. Then $R'_0\sub R_0$ and
$$\ti Y=\{x\in\fc;\{\al\in R_0;\al(x)=0\}\sub R'_0\}
=\{x\in\fc;\al(x)\ne 0\quad\frl\al\in R_0-R'_0\}.$$
This is non-empty, open in $\fc$; (a) is proved.

We prove (b). We can find a Levi subgroup $L$ of $Q$ such that 
$\fzg(z_s)\sub\ul$. Let $P\in\cp$ be such that $(z,P)\in\dy$. We have 
$z\in\up,\upi_P(z)\in\bocp+\fz_{\ubap}$ hence 
$z_s\in\up,\upi_P(z_s)\in\fz_{\ubap}$. Hence $z_s\in\upi_P\i(\fz_{\ubap})$,
that is $z_s$ is contained in a Cartan subalgebra of $\upi_P\i(\fz_{\ubap})$. 
Equivalently, $z_s\in c$ where $c$ is the centre of a Levi subalgebra $\up^1$ 
of $\up$. Since $c$ is abelian, from $z_s\in c$ we deduce $c\sub\fzg(x_s)$ 
hence $c\sub\ul$. Hence $\fz_{\ul}$ is contained in $\fzg(c)=\up^1$. Thus, 
$\fz_{\ul}\sub\up$ and (b) follows.

We prove (c). Let $(z,P)\in\dy$. Let $L,\fz_{\ul}$ be as in the proof of (b). 
We have $\fz_{\ul}\sub\up$. Since $z\in\fz_\fq(x_s)\sub\ul$, we have $z\in\ul$
hence $[\fz_{\ul},z]=0$. Then $[\upi_P(\fz_{\ul}),\upi_P(z)]=0$. By a known 
property of cuspidal local systems, the centralizer in $\ubap$ of an element in
$\bocp+\fz_{\ubap}$ (in particular, $\upi_P(z)$) has a unique Cartan 
subalgebra, namely $\fz_{\ubap}$. It follows that $\upi_P(\fz_{\ul})$ (the Lie
algebra of a torus) is contained in $\fz_{\ubap}$. Thus, 
$\fz_{\ul}\sub\upi_P\i(\fz_{\ubap})$. Hence $\fz_{\ul}$ is contained in a 
Cartan subalgebra of $\upi_P\i(\fz_{\ubap})$. Equivalently, $\fz_{\ul}\sub c'$
where $c'$ is the centre of a Levi subalgebra of $\up$. It follows that 
$\fzg(c')$ (a Levi subalgebra of $\up$) is contained in $\ul$ hence is 
contained in $\fq$. In particular, the $Q$-orbit of $P$ in $\cp$ is good. Thus,
$(z,P)\in\dy_w$ for some $w\in W_*$. This proves (c). 

We prove (d). Let $(z,P)\in\dy_w$. Let $L$ be as in (b). Then $Z^0_L\sub P$. By
2.1(b) (with $S$ replaced by $Z^0_L$) we have $(z,P^!)\in\dy_1$. Hence the 
morphism $f_w:\dy_w@>>>\dy_1$ as in (d) is well defined.

Assume that $(z,P),(z',P')$ in $\dy_w$ have the same image under $f_w$, that is
$z=z'$ and $P^!=P'{}^!$. Let $L$ be as in (b) (attached to $z=z'$). Then 
$Z^0_L\sub P,Z^0_L\sub P'$. Using 2.1(a) (with $S$ replaced by $Z^0_L$) we see
that $P^!=P'{}^!$ implies $P=P'$. Thus, $f_w$ is injective.

Now let $(z,P_1)\in\dy_1$. Let $P$ be the unique subgroup in $\boo(w)^S$ ($S$ 
as in 2.1) such that $P^!=P_1$. By 2.1(b) we have $(z,P)\in\dy_w$. We have
$f_w(z,P)=(z,P_1)$. Thus $f_w$ is surjective. We see that $f_w$ is bijective.
We omit the proof of the fact that $f_w\i$ is a morphism.

We prove (e). Let us first replace our cuspidal datum (see 1.4) by the cuspidal
datum $(\cb,\{0\},\bc)$ where $\cb$ is the variety of Borel subgroups of $G$.
Let $Y',W'_*,Y'_{w'}$ be the analogues of $\dy,W_*,\dy_w$ for this new cuspidal
datum ($Y$ is unchanged). Now $Y'_1=\{(z,B);z\in Y,B\in\cb,B\sub Q,z\in\un B\}$
and $pr_1:Y'_1@>>>Y$ is proper since $\{B;B\in\cb,B\sub Q\}$ is projective. 
Using the isomorphism $Y'_{w'}@>>>Y'_1$ (as in (d)), we deduce that
$pr_1:Y'_{w'}@>>>Y$ is proper for any $w'\in W'_*$. Hence in the cartesian
diagram
$$\CD Y'\tim_YY'_{w'}@>>>Y'_{w'}\\@VaVV  @V pr_1VV\\Y'@>pr_1>>Y\endCD$$
the map $a$ is proper. Hence the image under $a$ of
$\{(\xi,\xi')\in Y'\tim_YY'_{w'};\xi\in Y'_{w'}\}$ (a closed subset of
$Y'\tim_YY'_{w'}$) is closed in $Y'$. But this image is just $Y'_{w'}$. We see
that $Y'_{w'}$ is closed in $Y'$.

We now return to the cuspidal datum in 1.4. Let $m:\cb@>>>\cp$ be the morphism
given by $m(B)=P$ where $B\sub P$. This induces a map from the set of 
$Q$-orbits on $\cb$ to the set of $Q$-orbits on $\cp$, which can be viewed as a
map $\bar m:W'_*@>>>W_*$. Let $Y''=\{(z,P);z\in Y,P\in\cp,z\in\up\}$. Let
$m':Y'@>>>Y''$ be given by $m'(z,B)=(z,m(B))$. It is clear that $m'$ is a
proper morphism. Now $\dy$ is a subvariety of $Y''$. The restriction of $m'$ 
defines a proper morphism $m'':m'{}\i(\dy)@>>>\dy$. Since $Y'_{w'}$ is closed 
in $Y'$, we see that $Y'_{w'}\cap m'{}\i(\dy)$ is closed in $m'{}\i(\dy)$ (here
$w'\in W'_*$). Let $w\in W_*$. Since 
$\cup_{w'\in W'_*;\bar m(w')=w}Y'_{w'}\cap m'{}\i(\dy)$ is closed in
$m'{}\i(\dy)$ and $m''$ is proper, it follows that
$m''(\cup_{w'\in W'_*;\bar m(w')=w}Y'_{w'}\cap m'{}\i(\dy))$ is a closed subset
of $\dy$. It is clear that this subset is just $\dy_w$. Thus $\dy_w$ is closed
in $\dy$. Then $\cup_{w_1\in W_*}\dy_{w_1}$ is closed in $\dy$ hence its
complement $\dy_w$ is open in $\dy$. This proves (e). The lemma is proved.

\subhead 2.10\endsubhead
Let $\ti K$ be the direct image with compact support of $\dcl^*$ under 
$$pr_1:\{(z,P);z\in\fq/\fn;P\in\cp^*,z\in\up;\upi_P(z)\in\bocp+\uzp\}@>>>
\fq/\fn.$$
Let $\ti K'$ be the direct image with compact support of $\dcl^*$ under 
$$pr_1:\{(z,P);z\in\fq;P\in\cp^*,z\in\up;\upi_P(z)\in\bocp+\un{Z_{\bap}}\}@>>>
\dq$$
By \cite{\LI, 3.4(a)} applied to $\fq/\uuq$ instead of $\fg$, $\ti K$ is an
intersection cohomology complex (ICC) supported by 
$$\bar Y_0=\{z\in\fq/\uuq;\exists P\in\cp^*,z\in\up;\upi_P(z)\in\bocp+\uzp\}$$
(a closed irreducible subset of $\fq/\uuq$) with a canonical $W_J$-action.
Clearly, $\ti K'$ is the inverse image of $\ti K$ under the obvious vector 
bundle $\fq@>>>\fq/\uuq$. Hence $\ti K'$ is an ICC supported by
$$Y_0=\{z\in\fq;\exists P\in\cp^*,z\in\up;\upi_P(z)\in\bocp+\uzp\}$$
(a closed irreducible subset of $\fq$) with a canonical $W_J$-action. 

If $X$ is a subvariety of $\fq$, then $\ti K'|_{\dx_1}$ (a complex of sheaves
on $\dx_1$) has a $W_J$-action inherited from $\ti K'$ hence there is a natural
$W_J$-action on the hypercohomology
$$H^j_c(X,\ti K'|_{\dx_1})=H^j_c(\dx_1,\dcl^*).$$
From the definitions we see that, if $\io:\dx_1@>>>\dx$ is the imbedding, the 
induced map
$$\io^*:H^j_c(\dx,\dcl^*)@>>>H^j_c(\dx_1,\dcl^*)$$
is compatible with the $W_J$-actions ($W_J$ acts on $H^j_c(\dx,\dcl^*)$ as the
restriction of the $W$-action 1.9).

\subhead 2.11\endsubhead
Let 
$$Y_1=\{z\in\fq;\exists P\in\cp^*,z\in\up,\upi_P(z)\in\bocp+\uzp,
\fzg(z_s)\sub\up\}.$$
We have $Y_1\sub Y_0$. We show that 

(a) $Y_1$ is open dense in $Y_0$.
\nl
Clearly, $Y_1\ne\emp$. Since $X_1=\{z\in\fg;\exists P\in\cp,z\in\up,
\upi_P(z)\in\bocp+\uzp,\fzg(z_s)\sub\up\}$ is open in

$X_0=\{z\in\fg;\exists P\in\cp,z\in\up,\upi_P(z)\in cl(\bocp)+\uzp\}$
\nl
(see \cite{\LII, 7.1}) it follows that $X_1\cap Y_0$ is open in 
$X_0\cap Y_0=Y_0$. Hence to prove (a) it is enough to show that 
$X_1\cap Y_0=Y_1$. The inclusion $Y_1\sub X_1\cap Y_0$ is obvious. Conversely,
let $z\in X_1\cap Y_0$. Thus, $z\in\fq$ and there exist $P\in\cp,P'\in\cp^*$ 
such that

$z\in\up,z\in\up',\upi_P(z)\in\bocp+\uzp,\upi_{P'}(z)\in cl(\boc_{P'})+
\un{Z_{\bap'}},\fzg(z_s)\sub\up$.
\nl
We have $\upi_P(z_s)\in\uzp,\upi_{P'}(z_s)\in\un{Z_{\bap'}}$; hence there
exists a Levi subalgebra $\fl$ of $\fp$ and a Levi subalgebra $\fl'$ of $\fp'$
such that $z_s\in\fz_\fl,z_s\in\fz_{\fl'}$. Then $\fl\sub\fzg(z_s)$. This,
together with $\fzg(z_s)\sub\up$ implies $\fl=\fzg(z_s)$. We also have 
$\fl'\sub\fzg(z_s)$ and $\dim\fl'=\dim\fl=\dim\fzg(z_s)$ hence 
$\fl'=\fzg(z_s)$. Thus, $\fl=\fl'$. Let $c$ be the nilpotent orbit in $\fl$
that corresponds to $\boc_P$ under $\fl@>\sim>>\ubap$ and also to $\boc_{P'}$
under $\fl@>\sim>>\un{\bar P'}$. Now $\fl$ is a Levi subalgebra of 
$\fp\cap\fp'$ hence $z=z'+z''$ where $z'\in\fl,z''\in\un{U_{P\cap P'}}$ are 
uniquely determined. We also have $z=z_1+z_2+z_3$ with 
$z_1\in c,z_2\in\fz_\fl,z_3\in\uup$ and $z=z'_1+z'_2+z'_3$ with 
$z'_1\in cl(c),z'_2\in\fz_\fl,z'_3\in\un{U_{P'}}$. It follows that 
$z_3\in\uup\cap\up'\sub\un{U_{P\cap P'}}$ and
$z'_3\in\un{U_{P'}}\cap\up\sub\un{U_{P\cap P'}}$. By uniqueness of the
decomposition $z=z'+z''$ we then have $z_1+z_2=z'=z'_1+z'_2$. Taking nilpotent
parts we deduce $z_1=z'_1$. In particular, $z'_1\in c$. We see that 

$z\in\up',P\in\cp^*,\upi_{P'}(z)\in\boc_{P'}+\un{Z_{\bap'}},\fzg(z_s)\sub\up'$
\nl
hence $z\in Y_1$. This proves (a).

Since $Y_0$ is irreducible, from (a) we deduce that $Y_1$ is irreducible. Now 
$Y\cap Y_0$ is closed in $Y$. Since $Y_1\sub Y$, from (a) we deduce that $Y_1$
is open dense in $Y\cap Y_0$. In particular, $Y\cap Y_0$ is irreducible.

The image of $pr_1:\dy@>>>Y$ has image contained in $Y_0$. (Indeed, by 
2.9(c),(d), the image of $pr_1:\dy@>>>Y$ is contained in the image of 
$pr_1:\dy_1@>>>Y$.) Thus we have maps $pr_1:\dy@>>>Y\cap Y_0,
pr_1:\dy_w@>>>Y\cap Y_0$. Taking the direct image with compact support of
$\dcl^*$ under $pr_1:\dy@>>>Y\cap Y_0$ (resp. $pr_1:\dy_w@>>>Y\cap Y_0$), we
get a complex of sheaves $K'$ (resp. $K'_w$) on $Y\cap Y_0$. From 2.9(c),(e) we
have canonically $K'=\opl_{w\in W_*}K'_w$ in the derived category. Since 
$K'=K|_{Y\cap Y_0}$ where $K$ (as in 1.8) has a natural $W$-action, we see that
$K'$ has a natural $W$-action. On the other hand, $W$ acts on $W_*$ by 
$w:w_1\mto w*w_1$ where $w*w_1$ is the element of minimal length in $ww_1W_J$.

\proclaim{Lemma 2.12} For $w\in W$ and $w_1\in W_*$ we have 
$wK'_{w_1}=K'_{w*w_1}$.
\endproclaim
Clearly, $K'_1=\ti K'_{Y\cap Y_0}$. (Note that $Y\cap Y_0$ is an open dense 
subset of $Y_0$ since $Y$ is open in $\fq$ and $Y\cap Y_0\ne\emp$.) It follows
that $K'_1$ is an ICC supported by $Y\cap Y_0$. Using 2.9(d) we deduce that for
any $w_1\in W_*$, $K'_{w_1}$ is isomorphic to $K'_1$ in the derived category. 
Hence $K'_{w_1}$ is an ICC supported by $Y\cap Y_0$. Since 
$K'=\opl_{w_1\in W_*}K'_{w_1}$ we see that $K'$ is an ICC supported by 
$Y\cap Y_0$. It is therefore sufficient to check the equality 
$wK'_{w_1}=K'_{w*w_1}$ over the open dense subset $Y_1$ of $Y\cap Y_0$. Using
\cite{\LI, 3.2(a)} we see that 
$$(Y_1)\dot{}=\{(z,P);z\in Y_1;\fzg(z_s)\text{ is a Levi subalgebra of }\up\}.
$$
Then $W$ acts freely on $(Y_1)\dot{}$ by $w:(z,P)\mto(z,P')$ where $P'\in\cp$ 
is defined by the condition that $\fzg(z_s)$ is a Levi subalgebra of $\up'$ and
$(P,P')$ is in the good $G$-orbit on $\cp\tim\cp$ corresponding to $w$. The
decomposition $(Y_1)\dot{}=\sqcup_{w_1\in W_*}(Y_1)\dot{}_{w_1}$ clearly
satisfies

$w(Y_1)\dot{}_{w_1}=(Y_1)\dot{}_{w*w_1}$
\nl
for any $w\in W,w_1\in W_1$. Using the definition of $K'$ and of the $W$-action
on it we see that $wK'_{w_1}=K'_{w*w_1}$ holds over $Y_1$. The lemma is proved.

\proclaim{Lemma 2.13} Assume that $H^*_c(\cp_y,\dcl^*)\ne 0$. Let 
$\io:\dca''_1@>>>\dca''$ be the inclusion. The linear map
$$H^j_c(\dca'',\dcl^*)@>>>\bc[W]\ot_{\bc[W_J]}H^j_c(\dca''_1,\dcl^*)$$
defined by $\xi\mto\sum_{w\in W_*}w\ot\io^*(w\i\xi)$ is an isomorphism.
\endproclaim
Let $z\in\ca''$. We have $z=y+n+t$ where $n\in\fn$ and $t\in\ft$. 

By 2.8(g) we have $H^*_c(\cp^*_y,\dcl^*)\ne 0$ hence by \cite{\LI, 8.6} we have
$Eu(\cp^*_y)\ne 0$. The set $\{P\in\cp^*_y,t\in\up\}$ is the fixed point set of
a torus action on $\cp^*_y$ hence it has the same Euler characteristic as 
$\cp^*_y$; in particular, this set is non-empty.

Let $P\in\cp^*_y$ be such that $t\in\up$. Since $\fn\sub\up$ we have $z\in\up$.
Since $t\in\ft$, we have $[t,x]\in\uuq\sub\uup$ for all $x\in\fq$. In 
particular, $[t,x]\in\uup$ for all $x\in\up$ hence 
$$\upi_P(t)\in\uzp.\tag a$$ 
Thus, $\upi_P(z)=\upi_P(y)+\upi_P(t)\in\bocp+\uzp$ so that $z\in Y_0$. Now $z$
and $y+t$ have the same image in $\fq/\fn$ hence $z_s$ and $t$ have the same 
image in $\fq/\fn$ hence $z_s$ and $t$ are in the same $\Ad(U_Q)$-orbit. Since 
$t\in\ft_r$ we have $\fzg(t)\sub\fq$ hence $\fzg(z_s)\sub\fq$ so that $z\in Y$.
Thus, $z\in Y_0\cap Y$. We see that $\ca''\sub Y_0\cap Y$. Taking the direct 
image with compact support of $\dcl^*$ under $pr_1:\dca''@>>>\ca''$ (resp. 
$pr_1:\dca''_w@>>>\ca''$) we obtain a complex of sheaves $K''$ (resp. $K''_w$)
on $\ca''$. This is the same as $K'|_{\ca''}$ (resp. $K'_w|_{\ca''}$). Hence we
have canonically $K''=\opl_{w\in W_*}K''_w$ and from 2.12 we deduce that the 
$W$-action on $K''$ satisfies $wK''_{w_1}=K''_{w*w_1}$ for any $w\in W$ and any
$w_1\in W_*$. Hence we have $K''=\opl_{w\in W_*}wK''_1$. It follows that we 
have canonically 
$$H^j_c(\ca'',K'')=\opl_{w\in W_*}H^j_c(\ca'',K''_w)=
\opl_{w\in W_*}wH^j_c(\ca'',K''_1)$$ 
that is,
$$H^j_c(\dca'',\dcl^*)=\opl_{w\in W_*}H^j_c(\dca''_w,\dcl^*)=
\opl_{w\in W_*}wH^j_c(\dca''_1,\dcl^*).$$
The lemma follows.

\proclaim{Lemma 2.14} Let $\io':\dca'_1@>>>\dca'$ be the inclusion. The linear 
map
$$H^j_c(\dca',\dcl^*)@>>>\bc[W]\ot_{\bc[W_J]}H^j_c(\dca'_1,\dcl^*)$$
defined by $\xi\mto\sum_{w\in W_*}w\ot\io'{}^*(w\i\xi)$ is an isomorphism.
\endproclaim
If $H^*_c(\cp_y,\dcl^*)=0$ 
then the linear map above is $0@>>>0$, by 2.8(g) and
the result is obvious. Assume now that $H^*_c(\cp_y,\dcl^*)\ne 0$. It suffices
to prove the similar result where the ground field $\bc$ is replaced by an 
algebraic closure $\bar F_p$ of the finite field $F_p$ where $p$ is a large 
enough prime (local systems and cohomology will be $l$-adic, where $l$ is a
prime $\ne p$). We will assume (in this proof) that $G$ has a fixed $F_q$-split
rational structure ($F_q\sub\bar F_p$ has $q$ elements)
 with Frobenius map $F$,
that $y,\boc,Q,Q^1$ are $F$-stable and that we are given an isomorphism 
$F^*\cl@>\sim>>\cl$ which induces on any stalk at a point of $\boc(F_q)$ a map
of finite order. Moreover we assume that any nilpotent orbit in $\fg$ or $\fq$
is defined over $F_q$ and that any irreducible local system over such an orbit
is defined over $F_q$. Using the $l$-adic analogue of Lemma 2.13 and taking 
pure parts we obtain an isomorphism
$$H^j_c(\dca'',\dcl^*)_{pure}@>>>
\bbq[W]\ot_{\bbq[W_J]}H^j_c(\dca''_1,\dcl^*)_{pure}$$
where $H^j_c(?,?)_{pure}$ is the part of $H^j_c(?,?)$ where the Frobenius map 
acts with eigenvalues $\la$ such that any complex absolute value of $\la$ is 
$q^{j/2}$. It is then enough to show that
$$H^j_c(\dca'',\dcl^*)_{pure}=H^j_c(\dca',\dcl^*),\tag a$$
$$H^j_c(\dca''_w,\dcl^*)_{pure}=H^j_c(\dca'_w,\dcl^*)\tag b$$
for any $w\in W_*$. 

We prove (b). Using 2.5(b) (or rather its analogue over $\bar F_p$) we see that
it is enough to show that

$H^j_c(\dca''{}^1_1,\dcl^*)_{pure}=H^j_c(\dca'{}^1_1,\dcl^*).$
\nl
Using 2.5(a) (or rather its analogue over $\bar F_p$) we see that it is enough
to show that

$H^j_c((y+\ft_r)\tim\cp^*_y,\dcl^*)_{pure}=H^j_c((y+\ft)\tim\cp^*_y,\dcl^*).$
\nl
Using K\"unneth's theorem, we see that it is enough to show that
$$H^j_c(\cp^*_y,\dcl^*)_{pure}=H^j_c(\cp^*_y,\dcl^*),\tag c$$
$$H^j_c(\ft_r,\bbq)_{pure}=H^j_c(\ft,\bbq).\tag d$$
Since $\ft_r$ is the complement in $\ft$ of a finite set of hyperplanes, the
eigenvalues of $F$ on $H^j_c(\ft_r,\bbq)$ are easily seen to be of the form
$q^{j-\dim\ft}$ (see \cite{\LE}) and (d) follows.

Now (c) is a special case of
$$H^j_c(\cp_x,\dcl^*)_{pure}=H^j_c(\cp_x,\dcl^*),\tag e$$
(for any nilpotent element $x\in\fg(F_q)$) obtained by replacing $G$ by $Q^1$.
To prove (e) is the same as to prove that the eigenvalues of $F$ on the stalk 
at $x\in\fg(F_q)$ of the $j$-th cohomology sheaf of $K$ have complex absolute 
value $q^{j/2}$. Now the restriction of $K$ to the nilpotent variety of $\fg$ 
is of the form $\opl_{O,\ce}V_{O,\ce}\ot ICC(cl(O),\ce)$ where $O$ runs over 
the nilpotent orbits in $\fg$ and $\ce$ are irreducible local systems on $O$ 
that are linked to our fixed cuspidal datum by the generalized Springer 
correspondence \cite{\LIC, Sec. 6}; $V_{O,\ce}$ are certain multiplicity
spaces. It is then enough to prove that for any of these $ICC(cl(O),\ce)$, the 
eigenvalues of $F$ on the stalk at $x\in\fg(F_q)$ of the $j$-th cohomology 
sheaf have complex absolute value $q^{j/2}$ (this holds by \cite{\LCS, 24.6})
and that $F$ acts trivially on each multiplicity space $V_{O,\ce}$. These
multiplicity spaces can be viewed as multiplicity spaces of the various
irreducible representations of $W$ in the regular representation of $W$ hence
$F$ acts on them trivially. This proves (e) hence (c) and (b).

We prove (a). We can arrange the $Q$-orbits on $\cp$ in a sequence
$\co_1,\co_2,\dots,\co_n$ as in the proof of 2.8. Thus, 
$R_m=\co_1\cup\co_2\cup\dots\cup\co_m$ is closed in $\cp$ for any $m\in[1,n]$.
We set $R_0=\emp$. 

To prove (a) it is enough to show that
$$H^j_c(\dca''_{R_m},\dcl^*)_{pure}=H^j_c(\dca'_{R_m},\dcl^*),$$
for any $m\in[0,n]$. We argue by induction on $m$. For $m=0$ the result is
trivial. Assume now that $m\ge 1$ and that the result is known for $m-1$. We 
have a cohomology exact sequence 
$$0@>>>H^j_c(\dca''_{\co_m},\dcl^*)@>>>
H^j_c(\dca''_{R_m},\dcl^*)@>>>H^j_c(\dca''_{R_{m-1}},\dcl^*)@>>>0$$
(we use the fact that $\dca''_{\co_m}$ is empty if $\co_m$ is not good, see 
2.9(c), and is both open and closed in $\dca''$ if $\co_m$ is good, see 2.9(e).
Taking pure parts in this exact sequence gives again an exact sequence
$$0@>>>H^j_c(\dca''_{\co_m},\dcl^*)_{pure}@>>>H^j_c(\dca''_{R_m},\dcl^*)_{pure}
@>>>H^j_c(\dca''_{R_{m-1}},\dcl^*)_{pure}@>>>0.$$
This exact sequence together with the exact sequence 2.8(f) are the rows of the
commutative diagram
$$\CD 0@>>>H^j_c(\dca''_{\co_m},)_{pure}@>>>H^j_c(\dca''_{R_m},)_{pure}
@>>>H^j_c(\dca''_{R_{m-1}},)_{pure}@>>>0\\@. @VVV @VVV @VVV @.\\
0@>>>H^j_c(\dca'_{\co_m},)@>>>
H^j_c(\dca'_{R_m},)@>>>H^j_c(\dca'_{R_{m-1}},)@>>>0\endCD$$
where the vertical maps are induced by the obvious open imbeddings and
the symbol $\dcl^*$ is omitted in the notation. Now the left vertical map is an
isomorphism. (If $\co_m$ is good this follows from (b); if $\co_m$ is not good
this
follows from the fact that $H^j_c(\dca'_{\co_m},\dcl^*)=0$, see 2.7, and that 
$\dca''_{\co_m}=\emp$, see 2.9(c).) The right vertical map is an isomorphism by
the induction hypothesis. It follows automatically that the middle vertical map
is an isomorphism. This proves (a). The lemma is proved.

\proclaim{Lemma 2.15} Let $\io':\dca'_1@>>>\dca'$ be the inclusion. The 
$H^*_C$-linear map
$$\bc[W]\ot_{\bc[W_J]}H^C_*(\dca'_1,\dcl)@>>>H^C_*(\dca',\dcl)\tag a$$
defined by $w\ot\xi\mto w\io'_!(\xi)$ is an isomorphism.
\endproclaim
By \cite{\LI, 3.8} (which is applicable in view of 2.8(a)), the two sides of
(a) are finitely generated projective $H^*_C$-modules which after applying
$H^*_{\{1\}}\ot_{H^*_C}$ become the analogous objects with $C$ replaced by 
$\{1\}$. Now (a) is an isomorphism when $C$ is replaced by $\{1\}$ (we take the
transpose of the isomorphism in Lemma 2.14). We see that the lemma can be 
deduced from the following statement which is easily verified.

Let $\car$ be the polynomial algebra over $\bc$ in the indeterminates
$x_1,x_2,\dots,x_n$ graded by $\deg(x_i)=2$ for all $i$. Let $\ci$ be the ideal
of $\car$ generated by $x_1,x_2,\dots,x_n$. Let $M,M'$ be two $\bn$-graded free
$\car$-modules and let $f:M@>>>M'$ be an $\car$-linear map compatible with the
gradings such that $f$ induces an isomorphism $M/\ci M@>\sim>>M'/\ci M'$. Then
$f$ is an isomorphism.

The lemma is proved.

\medpagebreak

Let $\io:\dca_1@>>>\dca$ be the inclusion. Consider the $H^*_C$-linear map
$$\bh\ot_{\bh'}H^C_*(\cp_y^*,\dcl)@>>>H^C_*(\dca,\dcl)\tag b$$
given by the composition
$$\align&\bh\ot_{\bh'}H^C_*(\cp_y^*,\dcl)=
\bc[W]\ot_{\bc[W_J]}H^C_*(\cp_y^*,\dcl)@>1\ot p^*>>
\bc[W]\ot_{\bc[W_J]}H^C_*(\dca_1,\dcl)\\&@>a>>H^C_*(\dca,\dcl)\tag c
\endalign$$
where $a$ is given by $w\ot\xi\mto w\io_!(\xi)$ and $p:\dca_1@>>>\cb^*_y$ is 
the affine space bundle $(z,P)\mto P$.

\proclaim{Theorem 2.16 (Strong induction theorem)} The $W$-module structure
and $\bs$-module structure on $H^C_*(\dca,\dcl)$ define an $\bh$-module 
structure on $H^C_*(\dca,\dcl)$. Moreover, the map 2.15(b) is an $\bh$-linear
isomorphism.
\endproclaim
By the argument in 2.8(g), we have a natural $H^*_C$-linear imbedding 
$H^C_*(\cp_y,\dcl)@>>>H^C_*(\dca,\dcl)$ which becomes an isomorphism after the
scalars are extended to the quotient field of $H^*_C$. Since the $W$-module 
structure and $\bs$-module structure on $H^C_*(\cp_y,\dcl)$ are known to define
an $\bh$-module structure, the same must then hold for $H^C_*(\dca,\dcl)$ 
(which is projective over $H^*_C$). This proves the first assertion of the
theorem. 

To prove the second assertion we may assume by 2.8(g) that 
$H^*_c(\cp_y,\dcl^*)\ne 0$. The composition of $C\sub G\tim\bc^*$ with 
$G\tim\bc^*@>pr_2>>\bc^*$ is surjective, as we see from the Morozov-Jacobson 
theorem for $y\in\fq^1$. Hence the image under the induced homomorphism 
$H^*_{\bc^*}@>>>H^*_C$ of the generator $\bor$ is a non-zero element of $H^*_C$
denoted again by $\bor$. It is enough to show that $a$ in 2.15(c) is an
isomorphism. Recall from 2.8 that $b=\dim\ft$. We show that

(a) the map $(\io_1)_!:H^C_j(\dca_1,\dcl)@>>>H^C_{j+2b}(\dca'_1,\dcl)$ induced
by the inclusion $\io_1:\dca@>>>\dca'$ is injective and its image equals
$\bor^bH^C_j(\dca'_1,\dcl)$.
\nl
Recall that 
$$\dca'_1=\{(z,P)\in\dfg; z\in y+\fz_\fn(\ph(f_0))+\ft,P\sub Q\},$$
$$\dca_1=\{(z,P)\in\dfg; z\in y+\fz_\fn(\ph(f_0)),P\sub Q\}.$$
We have an isomorphism 
$$k:\dca_1\tim\ft@>\sim>>\dca'_1\tag b$$
given by 
$((y+n,P),t)\mto(y+n+t,P)$ where $n\in\fz_\fn(\ph(f_0)),t\in\ft$. (We use the 
fact that $\ft\sub\upi\i(\uzp)+\uup$ for any $P\in\cp^*$, see 2.13(a).) Note
that $k$ is $C$-equivariant where the action $C\tim\ft@>>>\ft$ is
$((g,\la),t)\mto\la^{-2}t$.

This implies, by \cite{\LI, 1.10(b)} that the image of $(\io_1)_!$ is 
$\bor^bH^C_j(\dca'_1,\dcl)$. Since $\bor^b\ne 0$ in $H^*_C$ and
$H^C_*(\dca'_1,\dcl)$ is projective over $H^*_C$, we have

$\dim\bor^bH^C_j(\dca'_1,\dcl)=\dim H^C_j(\dca'_1,\dcl)=
\dim H^C_j(\dca_1,\dcl)$.
\nl
(The second equality follows from (b) and \cite{\LI, 1.4(e)}.) Hence 
$(\io_1)_!$ must be an isomorphism onto $\bor^bH^C_j(\dca'_1,\dcl)$. This 
proves (a). 

As in 2.8(c), let $D=\dim\dca,D'=\dim\dca'$. Let
$D_1=\dim\dca_1,D'_1=\dim\dca'_1$ so that $D'_1=D_1+b$.

From 2.15 we have 
$H^C_j(\dca',\dcl)=\sum_{w\in W}wH^C_{j-2D'+2D'_1}(\dca'_1,\dcl)$ and
$\io_!:H^C_*(\dca,\dcl)@>>>H^C_*(\dca',\dcl)$ (induced by the inclusion
$\io:\dca@>>>\dca'$) is $W$-equivariant. Hence from (a) we can deduce that
$$\align&\bor^bH^C_j(\dca',\dcl)=
\sum_{w\in W}w\bor^bH^C_{j-2D'+2D'_1}(\dca'_1,\dcl)\sub
\sum_{w\in W}w(\io_1)_!H^C_{j-2D'+2D'_1}(\dca_1,\dcl)\\&
\sub\sum_{w\in W}w\io_!H^C_{j-2D'+2D'_1+2D-2D_1}(\dca,\dcl)\sub
\io_!H^C_{j-2D'+2D'_1+2D-2D_1}(\dca,\dcl).\endalign$$
Thus, 
$$\bor^bH^C_j(\dca',\dcl)\sub\io_!H^C_{j-2D'+2b+2D}(\dca,\dcl).$$
It follows that
$$\align&\dim\bor^bH^C_j(\dca',\dcl)=\dim H^C_j(\dca',\dcl)\le
\dim\io_!H^C_{j-2D'+2b+2D}(\dca,\dcl)\\&\le\dim H^C_{j-2D'+2b+2D}(\dca,\dcl).
\endalign$$
These inequalities must be equalities since 
$\dim H^C_j(\dca',\dcl)=\dim H^C_{j-2D'+2b+2D}(\dca,\dcl)$ (see 2.8(c)). It
follows that

(c) $\bor^bH^C_j(\dca',\dcl)=\io_!H^C_{j-2D'+2b+2D}(\dca,\dcl)$ and $\io_!$ is
an isomorphism onto its image.
\nl
From 2.15(a) we deduce
$$\bc[W]\ot_{\bc[W_J]}\bor^bH^C_j(\dca'_1,\dcl)@>\sim>>
\bor^bH^C_{j+2D'-2D'_1}(\dca',\dcl)$$
which by (a),(c) becomes
$$\bc[W]\ot_{\bc[W_J]}H^C_j(\dca_1,\dcl)@>\sim>>H^C_{j-2D_1+2D}(\dca,\dcl).$$
The theorem is proved.

\subhead 2.17\endsubhead
Let $V$ be a finite dimensional $\bc$-vector space with an algebraic action
of $C$. let $[V]\in H^*_C$ be the element corresponding (as in 
\cite{\LI, 1.11}) to the regular function $\un C@>>>\bc$, $\xi\mto\det(\xi,V)$.
Here $\xi:V@>>>V$ is given by the associated Lie algebra representation of 
$\un C$ on $V$. Now $E=\fz_\fn(\ph(f_0))$ is a $C$-module for the restriction 
of the $G\tim\bc^*$-action on $\fg$. Hence $[E]$ is a well defined element of 
$H^*_C$.

Note that $E=\ker(ad(f_0):\fn@>>>\fn)\cong\cok(ad(e_0):\fn@>>>\fn)$
canonically. (Indeed, we have
$\fn=\ker(ad(f_0):\fn@>>>\fn)\opl\Im(ad(e_0):fn@>>>\fn)$ since $\fn$ is an
$\fs\fl_2(\bc)$-module.)

\proclaim{Lemma 2.18} (a) The homomorphism 
$H^C_*(\cp^*_y,\dcl)@>>>H^C_*(\dca_1,\dcl)$ induced by the inclusion 
$\cp^*_y\sub\dca_1$ is injective.

(b) The homomorphism $H^C_*(\cp_y,\dcl)@>>>H^C_*(\dca,\dcl)$ induced by the
inclusion $\cp_y\sub\dca$ is injective.

(c)  We have $[E]H^C_*(\dca_1,\dcl)\sub H^C_*(\cp^*_y,\dcl)$.

(d)  We have $[E]H^C_*(\dca,\dcl)\sub H^C_*(\cp_y,\dcl)$.
\endproclaim
(b) has already been noted at the end of 2.8. An entirely similar proof yields
(a).

We prove (c). We have an isomorphism 
$$\fz_\fn(\ph(f_0))\tim\cp_y^*@>\sim>>\dca_1$$
given by $(n,P)\mto(y+n,P)$ (we use the fact that 
$\fz_\fn(\ph(f_0))\sub\fn\sub\uup$ 
for any $P\in\cp^*$). Under this isomorphism, the 
inclusion $\cp^*_y\sub\dca_1$ corresponds to the map 
$\cp^*_y@>>>\fz_\fn(\ph(f_0))\tim\cp_y^*$ given by $P\mto(0,P)$. Now the result
follows using \cite{\LI, 1.10(b)}.

We prove (d). Using 2.16, (c) and the fact that the homomorphism in (b) is
$W$-equivariant, we see that
$$\align&[E]H^C_*(\dca,\dcl)=[E]\sum_{w\in W}wH^C_*(\dca_1,\dcl)=
\sum_{w\in W}w[E]H^C_*(\dca_1,\dcl)\\&\sub
\sum_{w\in W}wH^C_*(\cp_y,\dcl)\sub H^C_*(\cp_y,\dcl).\endalign$$
The lemma is proved.

\subhead 2.19\endsubhead
From 2.18 we see that we have a natural isomorphism

$H^*_C[[E]\i]\ot_{H^*_C}H^C_*(\cp_y,\dcl)@>\sim>>
H^*_C[[E]\i]\ot_{H^*_C}H^C_*(\dca,\dcl)$.
\nl
We combine this with 2.16; theorem 1.17 follows.

\head 3. Proof of Theorems 1.15, 1.21, 1.22\endhead
\subhead 3.1\endsubhead
Until the end of 3.8, let $P,L,T,W(L),R,\cc$ be as in 1.6. Let 
$R_1=\{\al\in R;2\al\notin R\}$. Then $R_1$ is (reduced) root system in 
$\ut^*$ with Weyl group $W(L)$. For $i\in I$, the set 
$R_1\cap\{\al_i,2\al_i\}$
consists of a single element; we call it $\al'_i$. Then $\{\al'_i;i\in I\}$ is
a set of simple roots for $R_1$. We have $\fg^{\al'_i}\sub\uup$.

\subhead 3.2\endsubhead
Let $\ph_0:\fs\fl_2(\bc)@>>>\ul$ be a Lie algebra homomorphism such that
$\ph_0(e_0)\in\cc$. (Such $\ph_0$ exists by the Morozov-Jacobson theorem.) Let 
$Z=Z(\Im(\ph_0))^0$. (A connected reductive subgroup of $G$.) Then
$\uz=\fz(\Im(\ph_0))$. 

\proclaim{Lemma 3.3} (a) $T$ is a maximal torus of $Z$.

(b) Let $N_Z(T)$ be the normalizer of $T$ in $Z$ so that $N_Z(T)/T$ is the Weyl
group of $Z$. The obvious homomorphism $N_Z(T)/T@>>>N(T)/L=W(L)$ is an 
isomorphism.

(c) We have $\fg^0\cap\uz=\ul\cap\uz=\ut$. For any $\al\in R_1$ we have
$\dim(\fg^\al\cap\uz)=1$. Hence $R_1$ is the root system of $\uz$ with respect
to $\ut$. 

(d) $\up\cap\uz$ is a Borel subalgebra of $\uz$.

(e) The map $J\mto\un{P_J}\cap\uz$ is a bijection between $\{J;J\sub I\}$ and
the set of parabolic subalgebras of $\uz$ that contain $\up\cap\uz$.

(f) $\fz_{\uz}=\fz_{\fg}$.
\endproclaim
For (a) see \cite{\LI, 2.6(a)}. For (b) see \cite{\LPE, 11.7(b)}. For (c) see
\cite{\LI, 2.9}. For (d) see \cite{\LPE, 11.7(a)}. 

We prove (e). The map in (e) is well defined by (d). This is a map between two
finite sets of the same cardinal, $2^{\sha(I)}$. To show that it is bijective,
it is enough to show that it is injective. Let $J,J'$ be two subsets of $I$ 
that satisfy $\un{P_J}\cap\uz=\un{P_{J'}}\cap\uz$. We must show that $J=J'$. We
have $\un{P_J}\cap\un{P_{J'}}=\un{P_{J\cap J'}}$. Then 
$\un{P_{J\cap J'}}\cap\uz=\un{P_J}\cap\uz=\un{P_{J'}}\cap\uz$ and it is enough
to show that $J\cap J'=J$ and $J\cap J'=J'$. Thus we are reduced to the case 
where $J\sub J'$. Assume that $J\ne J'$. Let $i\in J'-J$. Then 
$\un{P_i}\sub\un{P_{J'}}$, $\un{P_i}\not\sub\un{P_J}$. Using 1.6(a) with 
$Q=P_J$ we see that $E\cap\un{P_J}=0$ where $E=\fg^{-\al_i}\opl\fg^{-2\al_i}$.
Note that $E\sub\un{P_{J'}}$. Let $E'=E\cap\uz$. We have 
$E'\cap(\un{P_J}\cap\uz)=0$, $E'\sub\un{P_{J'}}\cap\uz$. Since $\dim(E')=1$ (by
(c)), we deduce that $\un{P_J}\cap\uz\ne\un{P_{J'}}\cap\uz$, a contradiction. 
This proves (e).

We prove (f). Since $\fz_{\fg}\sub\fz_{\uz}$, it is enough to show that the two
centres have the same dimension. From (a) and (c) we see that 
$\dim\fz_{\uz}=\dim(\ut)-\sha(I)$. It is easy to see that 
$\dim\fz_{\fg}=\dim(\ut)-\sha(I)$. This proves (f). The lemma is proved.

\subhead 3.4\endsubhead
Let $V\in\ci$ and let $D$ be the unique $\up$-stable line in $V$. For any 
$v\in D,x\in\ut$ we have $xv=\xi_V(x)v$ where $\xi_V\in\ut^*$ corresponds under
the obvious isomorphism $\ut@>\sim>>\up/[\up,\up]=\fh$ to the vector of $\fh^*$
denoted in 1.19 again by $\xi_V$. Now $\{x\in\fg;xD\sub D\}=\un{P_K}$ for a 
well defined $K\sub I$. We then say that $V\in\ci_K$. We say that $V\in\ci_K^0$
if $V\in\ci_K$ and $\fz_{\fg}$ acts as $0$ on $V$.

Let $i\in I$. Then $\{\xi_V; V\in\ci^0_{I-\{i\}}\}=
\{\vp_i,2\vp_i,3\vp_i,\dots\}$ where $\vp_i\in\ut^*$ (or $\vp_i\in\fh^*$) is 
well defined; we have $\vp_i=\xi_{\La^i}$ where $\La^i\in\ci^0_{I-\{i\}}$ is 
well defined up to isomorphism. 

\subhead 3.5\endsubhead
Let $\ci'{}^0$ be the collection of all simple finite dimensional $\uz$-modules
on which $\fz_{\uz}$ acts as $0$. Let $V'\in\ci'{}^0$ and let $D'$ be the 
unique $(\up\cap\uz)$-stable line in $V$. There is a unique vector 
$\xi'_{V'}\in\ut^*$ such that $xv=\xi'_{V'}(x)v$ for any $v\in D',x\in\ut$.
Also $\{x\in\uz;xD'\sub D'\}=\un{P_K}\cap\uz$ for a well defined $K\sub I$ (see
3.3(e)). We then say that $V'\in\ci'{}^0_K$. 

Let $i\in I$. We have $\{\xi_{V'};V'\in\ci'{}^0_{I-\{i\}}\}=
\{\vp'_i,2\vp'_i,3\vp'_i,\dots\}$ where $\vp'_i\in\ut^*$ is well defined; we 
have $\vp'_i=\xi'_{\La'_i}$ where $\La'_i\in\ci'{}^0_{I-\{i\}}$ is well defined
up to isomorphism.

\proclaim{Lemma 3.6} Let $i\in I$. We have $\vp_i=n_i\vp'_i$ for some 
$n_i\in\bn-\{0\}$.
\endproclaim
Let $D$ be the unique $\up$-stable line in $\La^i$. Then 
$\un{P_{I-\{i\}}}=\{x\in\fg;xD\sub D\}$. We may regard $\La^i$ as a 
$\uz$-module by restriction. In this $\uz$-module, $\fz_{\uz}$ acts as $0$ (see
3.3(f)) and $D$ is stable under the Borel subalgebra $\up\cap\uz$ (see 3.3(d))
of $\uz$. Hence the $\uz$-submodule $V'$ generated by $D$ is simple. Clearly, 
$\{x\in\uz;xD\sub D\}=\un{P_{I-\{i\}}}\cap\uz$. Thus, 
$V'\in\ci'{}^0_{I-\{i\}}$. From the definition, we then have 
$\vp_i\in\{\vp'_i,2\vp'_i,3\vp'_i,\dots\}$. The lemma is proved.

\subhead 3.7\endsubhead
Let $X=\Hom(\bt,\bc^*)$ (homomorphisms of algebraic groups). The differential
gives an imbedding $d:X@>>>\fh^*$ whose image $\cx$ is a free abelian group 
such that $\bc\ot\cx@>\sim>>\fh^*$. Let
$\fh_\br=\{x\in\fh;\xi(x)\in\br\quad\frl\xi\in\cx\}$.

Under the obvious isomorphism $\ut@>\sim>>\up/[\up,\up]=\fh$, $\fh_\br$
corresponds to a subset $\ut_\br$ of $\ut$.

We shall need the following variant of a lemma of Langlands.

\proclaim{Lemma 3.8}Assume that $G$ is semisimple.

(a) For any $f\in\fh_\br$ there is a subset $J$ of $I$ and a decomposition 
$f={}^0f+{}^1f$ with ${}^0f,{}^1f\in\ut_\br$ such that
$$\al_i({}^0f)<0 \text{ if } i\in I-J,\quad \al_i({}^0f)=0 \text{ if }i\in J,$$
$$\vp_i({}^1f)\ge 0\text{ if }i\in J,\quad\vp_i({}^1f)=0 \text{ if }i\in I-J.$$
Moreover, $J,{}^0f,{}^1f$ are uniquely determined by $f$.

(b) If $f,f'\in\fh_\br$ satisfy $\vp_i(f)\le\vp_i(f')$ for any $i\in I$ then
$\vp_i({}^0f)\le\vp_i({}^0f')$ for any $i\in I$.
\endproclaim
Using the isomorphism $\ut@>\sim>>\fh$ in 3.7, we see that the statement above
is equivalent to the one where $\fh,\fh_\br$ are replaced by $\ut,\ut_\br$.
Moreover, if $\al_i,\vp_i$ are replaced by $\al'_i,\vp'_i$, then these
statements hold by Langlands' lemma \cite{\BW, IV, 6.11-6.13} applied to the 
root system $R_1$ of 
$\uz$ with respect to $\ut$. However this replacement does
not affect the statements since $\al'_i,\vp'_i$ differ from $\al_i,\vp_i$ only
by rational $>0$ factors (see 3.1, 3.6). The lemma is proved.

\subhead 3.9\endsubhead
In the remainder of this section (except in 3.42) we fix $y\in\fg_N$, $r\in\bc$
and $\si\in\fg_{ss}$ with $[\si,y]=2ry$ such that
$$Eu(\cp^\si_y)\ne 0\tag a$$
where
$$\cp^\si=\{P\in\cp;\si\in\up\},\cp^\si_y=\cp^\si\cap\cp_y$$ 
and $Eu$ denotes Euler characteristic. Condition (a) is equivalent to each of 
the following conditions.
$$Euler(\cp_y)\ne 0,\tag b$$
$$H_*(\cp_y,\dcl)\ne 0,\tag c$$
The equivalence of (a),(b) follows from the conservation of $Eu$ by passage to
the fixed point set of a torus action. The equivalence of (b),(c) follows from
the vanishing theorem \cite{\LI, 8.6}. 

If $r\ne 0$, conditions (a),(b),(c) are also equivalent to the condition
$$E_{y,\si,r}\ne 0.\tag d$$
The equivalence of (c),(d) follows from \cite{\LI, 7.2} (which is applicable in
view of \cite{\LI, 8.6}).

\proclaim{Lemma 3.10} Let $Q\in\pa$ and let $Q^1$ be a Levi subgroup of $Q$.
Assume that $y,\si$ are contained in $\uq^1$. Then 
$\{P\in\cp^\si_y;P\sub Q\}\ne\emp$.
\endproclaim
We may use 
the argument in the second paragraph of the proof of 2.13 (replacing
$t$ by $\si$).

\subhead 3.11\endsubhead
By a variant of the Morozov-Jacobson theorem (see \cite{\KL, 2.4(g)}) we can 
find $h,\ty$ in $\fg$ such that
$$[\si,\ty]=-2r\ty,[y,\ty]=h,[h,y]=2y,[h,\ty]=-2\ty.$$
Then $[\si,h]=0$. 

We now assume that $r\ne 0$. We fix $\ta:\bc@>>>\br$ as in 1.20.

Let $V$ be a finite dimensional $\fg$-module. We have 
$$V=\opl_{a\in\bc}V_a \text{ where } V_a=\{x\in V;\si x=ax\},$$
$$V=\opl_{n\in\bz}({}_nV) \text{ where } {}_nV=\{x\in V; hx=nx\}.$$
Since $[\si,h]=0$, the maps $v\mto\si v,v\mto hv$ from $V$ to $V$ commute hence
$$V=\opl_{n,a;n\in\bz,a\in\bc}({}_nV_a) \text{ where }{}_nV_a={}_nV\cap V_a.$$
We have
$$V=\opl_{b\in\br}{}^bV \text{ where }{}^bV=\opl_{n,a;\ta(a)/\ta(r)=n+b}
({}_nV_a).$$
These definitions can be applied in particular with $V$ replaced by $\fg$ with
the $\ad$ action of $\fg$. We have 
$$y\in{}_2\fg_{2r},h\in{}_0\fg_0,y'\in{}_{-2}\fg_{-2r},\si\in{}_0\fg_0.$$
From the definition we have
$$x\in{}_n\fg_a,v\in{}_{n'}V_{a'}\impl xv\in{}_{n+n'}V_{a+a'},$$
$$x\in{}^b\fg,v\in{}^{b'}V\impl xv\in{}^{b+b'}V.$$

\subhead 3.12\endsubhead
We define $Q\in\pa$ and a Levi subgroup $Q^1$ of $Q$ by
$$\uq=\opl_{n,a;\ta(a)/\ta(r)\le n}({}_n\fg_a)=\opl_{b;b\le 0}{}^b\fg,$$
$$\uq^1=\opl_{n,a;\ta(a)/\ta(r)=n}({}_n\fg_a)=\opl_{b;b=0}{}^b\fg.$$
Then
$$\fn=\opl_{n,a;\ta(a)/\ta(r)<n}({}_n\fg_a)=\opl_{b;b<0}{}^b\fg$$
is the nil-radical of $\uq$. Also, $y,h,\ty,\si$ are contained in $\uq^1$. Now
$\ad(\si),\ad(y)$ define endomorphisms of $\fn$ whose commutator is $2r\ad(y)$.
Hence $\ad(\si)$ maps 
$${}_y\fn=\cok(\ad(y):\fn@>>>\fn)$$ 
into itself. 

\proclaim{Lemma 3.13} $\ad(\si)-2r:{}_y\fn@>>>{}_y\fn$ is invertible.
\endproclaim
An equivalent statement is that any eigenvalue of $\ad(\si)-2r$ on 
$$\cok(\ad(y):\opl_{n,a;\ta(a)/\ta(r)<n}({}_n\fg_a)@>>>
\opl_{n,a;\ta(a)/\ta(r)<n}({}_n\fg_a)\tag a$$
is $\ne 0$. Now (a) is a quotient of 
$\opl_{n,a;n\le 0,\ta(a)/\ta(r)<n}({}_n\fg_a)$ which is itself a quotient of
$$\opl_{a;\ta(a)/\ta(r)<0}\fg_a\tag b$$
and it is enough to show that any eigenvalue $\la$ of $\ad(\si)-2r$ on (b) is 
$\ne 0$. We have $\la=a-2r$ for some $a\in\bc$ such that $\ta(a)/\ta(r)<0$. 
Then $\ta(\la)/\ta(r)=\ta(a-2r)/\ta(r)=\ta(a)/\ta(r)-2<-2$. In particular,
$\la\ne 0$.

\proclaim{Lemma 3.14} Let $z\in G$ be such that $\Ad(z)y=2cy,\Ad(z)\si=\si$ for
some $c\in\bc$. Then $z\in Q$.
\endproclaim
We must show that $\Ad(z)x\in\uq$ for any $x\in\uq$. We may assume that 
$x\in{}_n\fg_a$ where $\ta(a)/\ta(r)\le n$. We have 
$[\si,\Ad(z)x]=\Ad(z)[\si,x]=\Ad(z)(ax)=a\Ad(z)x$. Thus, $\Ad(z)x\in\ug_a$.
Since $\Ad(z)y=2cy$ for some $c$, $z$ belongs to the parabolic subgroup of $G$
corresponding to $\opl_{m\ge 0}({}_m\fg)$. Hence
$\Ad(z)(x)\in\Ad(z)({}_n\ug)\sub\opl_{m;m\ge n}({}_m\fg)$. We see that
$$\Ad(z)(x)\in\opl_{m;m\ge n}({}_m\fg_a)\sub\opl_{m;\ta(a)/\ta(r)\le m}
({}_m\fg_a)\sub\uq.$$
The lemma is proved.

\proclaim{Lemma 3.15} (a) $Q$ is independent of the choice of $h,\ty$.

(b) We have $M(y,\si)\sub Q\tim\bc^*$.
\endproclaim
We prove (a). Any other choice of $h,\ty$ is of the form $h',\ty'$ where 
$h'=\Ad(z)h$, $\ty'=\Ad(z)\ty$ for some $z\in G$ such that 
$\Ad(z)y=y,\Ad(z)\si=\si$. (See \cite{\KL, 2.4(h)}.) Let $Q'$ be attached to 
$h',\ty'$ in the same way as $Q$ is attached to $h,\ty$. Then $Q'=zQz\i$. By 
3.14 we have $z\in Q$. Hence $Q'=Q$.

We prove (b). Let $(g,\la)\in M(y,\si)$. We can find an element $g_1$ in the 
one parameter subgroup of $G$ corresponding to $h$ such that 
$\Ad(g_1)y=\la^2y$. Since $h\in\uq^1,[\si,h]=0$, we have 
$g_1\in Q^1\cap Z(\si)$. Replacing $(g,\la)$ by $(gg_1\i,1)$ we see that we are
reduced to the case where $\la=1$ and $g\in Z(y)\cap Z(\si)$.

Since $g\in Z(\si)$ we have $\Ad(g)(\fg_a)=\fg_a$ for all $a$. Since
$g\in Z(y)$, we have $\Ad(g)({}_n\fg)\sub\opl_{n';n'\ge n}({}_{n'}\fg)$. Hence
$\Ad(g)({}_n\fg_a)\sub\opl_{n';n'\ge n}({}_{n'}\fg_a)$ for any $n,a$. Using
this and the definition of $\uq$ (see 3.12) we see that $\Ad(g)(\uq)\sub\uq$.
Hence $g\in Q$. The lemma is proved.

\proclaim{Lemma 3.16} Let $\si',y',h',\ty'$ be another quadruple like 
$\si,y,h,\ty$. Define $Q'$ in terms of $\si',y',h',\ty'$ in the same way as
$Q$ was defined in terms of $\si,y,h,\ty$. Assume that $Q=Q'=G$. Assume that
there exist $P\in\cp^\si_y$ and $P'\in\cp^{\si'}_{y'}$ such that the image of
$\si$ in $\up/[\up,\up]=\fh$ coincides with the image of $\si'$ in 
$\up'/[\up',\up']=\fh$. Then there exists $g\in G$ such that $\Ad(g)$ carries 
$(\si,y,h,\ty)$ to $(\si',y',h',\ty')$.
\endproclaim
Replacing $(\si',y',h',\ty')$ by a $G$-conjugate we may assume that $P'=P$ and

$\pi_P(\si)=\si_1,\pi_P(\si_1)=\si'_1$, $\pi_P(y')=\pi_P(y)=y_1$,

$\pi_P(h')=\pi_P(h)=h_1$, $\pi_P(\ty')=\pi_P(\ty)=\ty_1$,
\nl
where $y_1\in\bocp$ and

$[y_1,\ty_1]=h_1,[h_1,y_1]=2y_1,[h_1,\ty_1]=-2\ty_1$,

$[\si_1,y_1]=2ry_1,[\si_1,\ty_1]=-2r\ty_1$,

$[\si'_1,y_1]=2ry_1,[\si'_1,\ty_1]=-2r\ty_1$.
\nl
Moreover $\si_1$ and $\si'_1$ have the same image in $\ubap/[\ubap,\ubap]$. 
Hence $x=\si_1-\si'_1\in[\ubap,\ubap]$. We have $[x,y_1]=0,[x,\ty_1]=0$ hence
$[x,h_1]=0$. Since $y_1$ is a distinguished nilpotent element of $\ubap$, the
centralizer in $[\ubap,\ubap]$ of $y_1,h_1,\ty_1$ is $0$. Thus, $x=0$ so that
$\si_1=\si'_1$. Let $L$ be a Levi subgroup of $P$. Since $\si,\si'\in\up$, 
there exist $g,g'\in P$ such that 
$$\Ad(g)\si\in\ul,\pi_P(\si-\Ad(g)\si)=0,\Ad(g')\si'\in\ul, 
\pi_P(\si'-\Ad(g')\si')=0.$$ 
Hence $\pi_P(\Ad(g)\si)=\si_1=\si'_1=\pi_P(\Ad(g')\si')$. Since the restriction
of $\pi_P$ to $\ul$ is injective, it follows that $\Ad(g)\si=\Ad(g')\si'$. 
Thus, $\si,\si'$ are conjugate in $G$.

Replacing $(\si',y',h',\ty')$ by a $G$-conjugate we may assume that $\si'=\si$.
We show that 

(a) {\it $y$ belongs to the (unique) open orbit of $Z(\si)$ on $\fg_{2r}$.}
\nl
Let $G'$ be the connected reductive algebraic subgroup of $G$ such that
$$\ug'=\opl_{m\in\bz}\fg_{2mr}.$$
Note that $y,h,\ty$ are contained in $\ug'$. Let ${}_n\ug'={}_n\fg\cap\ug'$.
Since in our case, $\uq^1=\fg$, any eigenvalue $b$ of $\ad(\si-rh):\fg@>>>\fg$
satisfies $\ta(b)=0$. Hence ${}_n\fg_a\ne 0\impl\ta(a)/\ta(r)=n$. Hence 
$\fg_0={}_0\ug'$ and $\fg_{2r}={}_2\ug'$. It follows that $Z(\si)$ is equal to
the centralizer $C(h)$ of $h$ in $G'$. We are reduced to the following known 
statement about $\fs\fl_2$-triples in $\ug'$: $y$ belongs to the open orbit of
$C(h)$ on ${}_2\ug'$. Thus, (a) holds.

Similarly, $y'$ belongs to the (unique) open orbit of $Z(\si')$ on 
$\{x\in\fg;[\si',x]=2rx\}$. Since $\si=\si'$, we see that both $y$ and $y'$ 
belong to the unique open orbit of $Z(\si)$ on $\fg_{2r}$. In particular,
$y,y'$ are conjugate under an element in $Z(\si)$. 

Replacing $(\si,y',h',\ty')$ by a $Z(\si)$-conjugate we may therefore assume
that $y=y'$. As in the proof of 3.15, we can find $z\in G$ such that 
$\Ad(z)y=y,\Ad(z)\si=\si,h'=\Ad(z)h,\ty'=\Ad(z)\ti y$. Replacing 
$(\si,y,h',\ty')$ by its $\Ad(z\i)$-conjugate we may therefore assume that 
$(\si',y',h',\ty')=(\si,y,h,\ty)$. The lemma is proved.

\proclaim{Lemma 3.17} We have $\{P\in\cp^\si_y;P\sub Q\}\ne\emp$. In 
particular, $Q\in\cp_K$ for a well defined $K\sub I$.
\endproclaim
This follows from 3.10.

\subhead 3.18\endsubhead
Let $V\in\ci_J$ where $J\sub I$. For any $P\in\cp$ let $D_P$ be the unique 
$\up$-stable line in $V$. If $P\in\cp^\si$, we have necessarily $D_P\sub V_a$
for a unique $a\in\bc$; we set $\nu_V(P)=a$. Let
$$b_V=\min(b\in\br;{}^bV\ne 0).$$
The following result is closely related to \cite{\LSQ, 2.8, 2.9}, 
\cite{\KL, 7.3}.

\proclaim{Lemma 3.19} We preserve the setup of 3.18.

(a) If $P\in\cp^\si,y\in\up$, then $\ta(\nu_V(P))/\ta(r)\ge b_V$.

(b) If $K=J,P\in\cp^\si,y\in\up,P\not\sub Q$, then $\ta(\nu_V(P))/\ta(r)>b_V$.

(c) If $P\in\cp,P\sub Q$, then $D_P\sub{}^{b_V}V$.

(d) If $K=J,P\in\cp^\si,P\sub Q$, then $\ta(\nu_V(P))/\ta(r)=b_V$.
\endproclaim
From 3.11 we have
$$\fn{}^bV\sub\opl_{b';b'<b}{}^{b'}V,\uq{}^bV\sub\opl_{b';b'\le b}{}^{b'}V.$$
Hence $\fn{}^{b_V}V=0,\uq{}^{b_V}V\sub{}^{b_V}V$.

Let $P$ be as in (a) and let $v\in D_P-\{0\}$. We have $v\in V_{\nu_V(P)}$. We
write $v=\sum_m({}_mv)$ where ${}_mv\in{}_mV_{\nu_V(P)}$. Since $v\ne 0$, there
exists $n$ such that ${}_nv\ne 0$. Since $y$ is nilpotent in $\up$, we have
$yv=0$ hence $\sum_my({}_mv)=0$. Since $y\in{}_2\fg_{2r}$, we have
$$y({}_mv)\in{}_{m+2}V_{\nu_V(P)+2r}.$$ 
Since $\sum_my({}_mv)=0$ and the sum $\sum_m({}_{m+2}\fg_{\nu_V(P)+2r})$ is
direct, we have $y({}_mv)=0$ for all $m$. In particular, $y({}_nv)=0$. From 
$y({}_nv)=0$ and ${}_nv\ne 0$ we see, using the representation theory of 
$\fs\fl_2$ that $n\ge 0$. Since ${}_nv$ is a non-zero vector of 
${}_nV_{\nu_V(P)}\sub{}^{\ta(\nu_V(P))/\ta(r)-n}V$ and ${}^{b'}V=0$ unless 
$b'\ge b_V$, we see that $\ta(\nu_V(P))/\ta(r)-n\ge b_V$. Since $n\ge 0$, we 
must have $\ta(\nu_V(P))/\ta(r))\ge b_V$. This proves (a).

In the setup of (b), assume that $\ta(\nu_V(P))/\ta(r)\not>b_V$. Then, by (a),
we have $\ta(\nu_V(P))/\ta(r)=b_V$. Also, in the proof of (a) we must have
${}_nv\ne 0\impl n=0$ so that $v={}_0v\in{}_0V_{\nu_V(P)}$ and 
$v\in{}^{\ta(\nu_V(P))/\ta(r)}V={}^{b_V}V$. Thus, $D_P\sub{}^{b_V}V$. By 
assumption, $V$ contains a line $D$ such that $\{x\in\fg;xD\sub D\}=\uq$. This
implies that $\{v\in V;\fn v=0\}=D$. (See the argument in 1.19.) Since 
$\fn{}^{b_V}V=0$, it follows that ${}^{b_V}V\sub D$. Since $D_P\sub{}^{b_V}V$,
we have $D_P\sub D$ hence $D_P=D$. Since $D_P$ is $\up$-stable, we see that $D$
is $\up$-stable hence by the definition of $D$ we have $\up\sub\uq$ so that
$P\sub Q$. This proves (b).

Next, assume that $P$ is as in (c). Since $\uq{}^{b_V}V\sub{}^{b_V}V$, we have
$\up{}^{b_V}V\sub{}^{b_V}V$. Let $\fb$ be a Borel subalgebra of $\up$. Then 
$\fb{}^{b_V}V\sub{}^{b_V}V$ and, by Lie's theorem, there exists an $\fb$-stable
line $L$ in ${}^{b_V}V$. This is necessarily the unique $\fb$-stable line in 
$V$. Since $D_P$ is $\fb$-stable we must have $L=D_P$ hence $D_P\sub{}^{b_V}V$.
This proves (c).

In the setup of (d), $V$ contains a line $D$ such that 
$\{x\in\ug;xD\sub D\}=\uq$. Now $D$ is $\up$-stable, hence $D_P=D$. Since $h$ 
is contained in the derived subalgebra of $\uq$, it acts as zero on the 
$\uq$-stable line $D$. Hence $D\sub {}_0V$. We have $D_P\sub V_{\nu_V(P)}$ 
hence $D\sub V_{\nu_V(P)}$. As in the proof of (b) we have ${}^{b_V}V\sub D$; 
this must be an equality since $\dim{}^{b_V}V\ge 1,\dim D=1$. From
$D\sub{}_0V,D\sub V_{\nu_V(P)},D={}^{b_V}V$ we deduce 
$\ta(\nu_V(P))/\ta(r)=b_V$. This proves (d). The lemma is proved.

\proclaim{Lemma 3.20} In the setup of 3.18, assume that $Q=Q^1=G$. Let
$P'\in\cp_{J'}$ where $J'\sub J$ and let $L'$ be a Levi subgroup of $P'$. 
Assume that $\si,y,h,\ty$ are contained in $\ul'$. If $P\in\cp^\si,P\sub P'$,
then $\ta(\nu_V(P))/\ta(r)=b_V$.
\endproclaim
Clearly, $\uq^1{}^bV\sub{}^bV$ for any $b$. Since $\uq^1=\fg$, we see that
${}^bV$ is a $\fg$-submodule of $V$. Since $V$ is simple, we have $V={}^bV$ for
some $b$. Since ${}^{b_V}V\ne 0$, we have $V={}^{b_V}V$. Since $P'\in\cp_{J'}$,
there exists a $\up'$-stable line $D$ in $V$. From our assumptions, we have 
$h\in[\up',\up']$. Hence $hD=0$ so that $D\sub{}_0V$. Since $\si\in\up'$, we
have $\si D\sub D$ hence $D\sub V_a$ for some $a$. Thus, 
$D\sub{}_0V_a\sub{}^{\ta(a)/\ta(r)}V$. Thus, ${}^{\ta(a)/\ta(r)}V\ne 0$. Now 
${}^bV=0$ unless $b=b_V$. Hence $\ta(a)/\ta(r)=b_V$. Since $P\sub P'$, we have
$\up D\sub D$. Hence $D_P=D$ and $a=\nu_V(P)$. The lemma is proved.

\proclaim{Lemma 3.21} $\{P\in\cp^\si_y;P\sub Q\}$ is open and closed in 
$\cp^\si_y$.
\endproclaim
(Compare \cite{\KL, 7.4}.) We can find $V\in\ci_J$ with $J=K$. In terms of this
$V$ we define $\nu_V:\cp^\si@>>>\bc$ and $b_V$ as in 3.18. Since $\cp^\si$ is 
compact, $\nu_V:\cp^\si@>>>\bc$ is constant on any connected component of 
$\cp^\si$, hence it is locally constant. Hence its restriction 
$\nu_V:\cp^\si_y@>>>\bc$ is locally constant. Hence 
$\{P\in\cp^\si_y;\ta(\nu_V(P))/\ta(r)=b_V\}$ is open and closed in $\cp^\si_y$.
By 3.19(b),(d), we have 
$$\{P\in\cp^\si_y;\ta(\nu_V(P))/\ta(r)=b_V\}=\{P\in\cp^\si_y;P\sub Q\}.$$
The lemma follows.

\subhead 3.22\endsubhead
Let $\ps:\cp^\si@>>>\fh$ be the morphism whose value at $P$ is the image of 
$\si\in\up$ in $\up/[\up,\up]=\fh$. This must be locally constant since
$\cp^\si$ is compact and $\fh$ is affine. From the definitions, we have
$$\xi_V(\ps(P))=\nu_V(P)$$
for any $V\in\ci,P\in\cp^\si$. 

\proclaim{Lemma 3.23} If $V\in\ci_J$, $J=K$, 
$P\in\cp^\si,y\in\up,P\not\sub Q$ and $P'\in\cp^\si,P'\sub Q$, then
$$\ta(\xi_V(\ps(P)))/\ta(r)>\ta(\xi_V(\ps(P')))/\ta(r).$$
\endproclaim
In view of 3.22, an equivalent statement is 
$\ta(\nu_V(P))/\ta(r)>\ta(\nu_V(P'))/\ta(r)$ and this follows from 3.19(b),(d).

\subhead 3.24\endsubhead
Since $[\si,h]=0$, we have $\si-rh\in\fg_{ss}$. Let 
$\cp^{\si-rh}=\{P\in\cp;\si-rh\in\up\}$.

\proclaim{Lemma 3.25} Let $A=\{P\in\cp^{\si-rh};P\sub Q\}$. 

(a) We have $A\ne\emp$. 

(b) For $P\in A$, let $\btt_P$ be the image of $\si-rh$ in $\up/[\up,\up]=\fh$.
Let $V\in\ci$. We have $\ta(\xi_V(\btt_P))/\ta(r)=b_V$.

(c) $\btt_P\in\fh$ is independent of the choice of $P\in A$. We denote it by 
$\btt$.
\endproclaim
We prove (a). We have $\si-rh\in\uq$. Hence there exists a Borel subgroup $B$  
of $Q$ such that $\si-rh\in\un B$. By 3.17, we have
$\{P\in\cp;P\sub Q\}\ne\emp$. This is a conjugacy class of parabolic subgroups
of $Q$ hence at least one of its members contains $B$. This proves (a).

We prove (b). Let $D_P$ be as in 3.18. We have $V=\opl_{c\in\bc}{}^{(c)}V$ 
where 
$${}^{(c)}V=\opl_{n,a;a-rn=c}({}_nV_a)=\{x\in V;(\si-rh)x=cx\}.$$
We have ${}^bV=\opl_{c;\ta(c)/\ta(r)=b}{}^{(c)}V$. Define 
$\nu':\cp^{\si-rh}@>>>\bc$ by $\nu'(P)=c$ where $D_P\sub{}^{(c)}V$. For 
$P\in\cp^{\si-rh}$ we have $\xi_V(\btt_P)=\nu'(P)$. By 3.19(c), for 
$P\in\cp,P\sub Q$ we have $D_P\sub{}^{b_V}V$, that is, 
$D_P\sub\opl_{c;\ta(c)/\ta(r)=b_V}{}^{(c)}V$. Thus if $P\in A$, then 
$\nu'(P)=c$ where $\ta(c)/\ta(r)=b_V$. Hence $\xi_V(\btt_P)=c$ where 
$\ta(c)/\ta(r)=b_V$. Hence $\ta(\xi_V(\btt_P))=b_V\ta(r)$. This proves (b).

We prove (c). Let $P',P''$ be two members of $A$. Let $V\in\ci$. By (b), we
have $\ta(\xi_V(\btt_{P'}))=\ta(\xi_V(\btt_{P''}))$. Since this holds for any
$\ta$, we have $\xi_V(\btt_{P'})=\xi_V(\btt_{P''})$. Since $\xi_V$ with 
$V\in\ci$ span $\fh^*$, it follows that $\btt_{P'}=\btt_{P''}$. The lemma is 
proved.

\proclaim{Lemma 3.26} Let $P'\in\cp^\si$ be such that $y\in\up'$, $P'\sub Q$. 
(Such $P'$ exists by 3.17).

(a) If $i\in I-K$ then $-\ta(\al_i(\btt))/\ta(r)>0$.

(b) If $i\in K$ then $-\ta(\al_i(\btt))/\ta(r)=0$.

(c) Let $V\in\ci$. Then $\ta(\xi_V(\ps(P')-\btt))/\ta(r)\ge 0$.

(d) Let $V\in\ci_J$ where $K\sub J\sub I$. Then 
$\ta(\xi_V(\ps(P')-\btt))/\ta(r)=0$.
\endproclaim
Pick $P\in A$ (see 3.25(a)). Since $\si-rh$ is a semisimple element of
$\uq^1\cap\up$ and $P\sub Q$, we can find a Levi subgroup $L$ of $P$ such that
$L\sub Q^1$ and $\si-rh\in\ul$. Let $T=Z^0_L$. Under the obvious isomorphism
$\ut@>\sim>>\up/[\up,\up]=\fh$, $\btt\in\fh$ corresponds to an element 
$t\in\ut$ such that $t=\si-rh+x$ where $x\in[\up,\up]$. Since 
$\si-rh\in\ul,t\in\ul$, we have $x\in[\up,\up]\cap\ul=[\ul,\ul]$. Let $i\in I$.
Now $\fg^{-\al_i}$ (defined as in 1.6 in terms of our $P,L$) is an $\ul$-module
hence $\tr(x,\fg^{-\al_i})=0$ (since $x\in[\ul,\ul]$). Let $z_1,z_2,\dots,z_k$
be the eigenvalues of $\ad(\si-rh)$ on $\fg^{-\al_i}$. 

Assume first that $i\in K$. Then $\un{L_i}\sub\uq^1$ (notation of 1.6) hence 
$\fg^{-\al_i}\sub\uq^1$. Using this and the definition of $\uq^1$, we have 
$\ta(z_j)=0$ for all $j\in[1,k]$. Then $\tr(\si-rh,\fg^{-\al_i})=z_1+\dots+z_k$
and $\ta(\tr(\si-rh,\fg^{-\al_i}))=\ta(z_1)+\dots+\ta(z_k)=0+\dots+0=0$. Now 
$t$ acts on $\fg^{-\al_i}$ as $-\al_i(t)$ times the identity hence
$\tr(t,\fg^{-\al_i}))=-k\al_i(t)$. Since $t=\si-rh+x$, we have 
$-\ta(k\al_i(t))=0+0=0$ and $-\ta(\al_i(t))=0$. 

Assume next that $i\in I-K$. Then $\un{P_i}\not\sub\uq$ (notation of 1.6). 
Using 1.6(a) we see that $(\fg^{-\al_i}\opl\fg^{-2\al_i})\cap\uq=0$. Hence 
$\fg^{-\al_i}\cap\uq=0$. This implies that $\ta(z_j)>0$ for all $j\in[1,k]$.
Then $\tr(\si-rh,\fg^{-\al_i})=z_1+\dots+z_k$ and
$\ta(\tr(\si-rh,\fg^{-\al_i}))=\ta(z_1)+\dots+\ta(z_k)>0$. Now 
$\tr(t,\fg^{-\al_i}))=-k\al_i(t)$. Since $t=\si-rh+x$ we have 
$-\ta(k\al_i(t))>0$ and $-\ta(\al_i(t))>0$. This proves (a) and (b).

We prove (c). Using 3.25(b) and 3.22 we see that we must prove that
\linebreak $\ta(\nu_V(P'))/\ta(r)\ge b_V$. This follows from 3.19(a).

We prove (d). The subgroup generated by the $\xi_V$ with $V\in\ci_K$ contains 
any $\xi_V$ with $V\in\cp_J$ where $K\sub J\sub I$. Hence it suffices to prove
(d) for $V\in\ci_K$. Using 3.25(b) and 3.22 we see that we must prove that
$\ta(\nu_V(P'))/\ta(r)=b_V$ when $V\in\ci_K$. This follows from 3.19(d). The 
lemma is proved.

\subhead 3.27\endsubhead
For $x\in\fh$ define ${}^\ta x\in\fh_\br$ (see 3.7) by 
$\ga({}^\ta x)=\ta(\ga(x))/\ta(r)$ for all $\ga\in\cx$. Then $x\mto{}^\ta x$ is
a group homomorphism $\fh@>>>\fh_\br$. We can now reformulate Lemma 3.26 as 
follows.

\proclaim{Lemma 3.28} Let $P'\in\cp^\si$ be such that $y\in\up'$, $P'\sub Q$. 

(a) If $i\in I-K$ then $-\al_i({}^\ta\btt)>0$.

(b) If $i\in K$ then $-\al_i({}^\ta\btt)=0$.

(c) Let $V\in\ci$. Then $\xi_V({}^\ta\ps(P')-{}^\ta\btt)\ge 0$.

(d) Let $V\in\ci_J$ where $K\sub J\sub I$. Then
$\xi_V({}^\ta\ps(P')-{}^\ta\btt)=0$.
\endproclaim

\proclaim{Lemma 3.29} Assume that $G$ is semisimple. Let $P'\in\cp^\si$ be such
that $y\in\up'$, $P'\sub Q$. Then ${}^0({}^\ta\ps(P'))={}^\ta\btt\in\fh_\br$ 
(notation of 3.8(a)).
\endproclaim
This follows immediately from 3.8(a) and 3.28.

\proclaim{Lemma 3.30} Let $\si',y',h',\ty'$ be another quadruple like 
$\si,y,h,\ty$. Define \linebreak $Q',Q'{}^1,\btt'\in\fh$ in terms of 
$\si',y',h',\ty'$ in the same way that $Q,Q^1,\btt\in\fh$ were defined in terms
of $\si,y,h,\ty$. Define $K'\sub I$ by $Q'\in\cp_{K'}$. Assume that there exist
$P_1,P_2\in\cp^\si_y$ with $P_1\sub Q$ and $P'_1,P'_2\in\cp^{\si'}_{y'}$ with
$P'_1\sub Q'$, such that 

(i) the image $\et$ of $\si$ in $\up_1/[\up_1,\up_1]=\fh$ coincides with the 
image of $\si'$ in $\up'_2/[\up'_2,\up'_2]=\fh$ and

(ii) the image $\et'$ of $\si'$ in $\up'_1/[\up'_1,\up'_1]=\fh$ coincides with
the image of $\si$ in $\up_2/[\up_2,\up_2]=\fh$.
\nl
Then there exists $g\in G$ such that $\Ad(g)$ carries $(\si,y,h,\ty)$ to 
$(\si',y',h',\ty')$.
\endproclaim
The general case reduces easily to the case where $G$ is semisimple. We now 
assume that $G$ is semisimple. Applying 3.29 twice (once for $\si,y,h,\ty,P_1$
and once for $\si',y',h',\ty',P'_1$) we see that 
$${}^0({}^\ta\et)={}^\ta\btt, {}^0({}^\ta\et')={}^\ta\btt'.$$ 
(Notation of 3.8(a).) Applying 3.19(a) (reformulated with the aid of 3.22 and 
3.25(b)) to $\si,y,h,\ty,P_2$ (where $P_2$ is not necessarily contained in $Q$)
we see that   
$$\xi_V({}^\ta\et'-{}^\ta\btt)\ge 0\tag a$$
for any $V\in\ci$. Using (a) and 3.8(b) with $f={}^\ta\btt',f'={}^\ta\et$ we 
deduce that $\xi_{\La_i}({}^0({}^\ta\et'))\ge\xi_{\La_i}({}^0({}^\ta\btt))$ for
any $i\in I$. Hence $\xi_{\La_i}({}^\ta\btt')\ge\xi_{\La_i}({}^\ta\btt)$ for 
any $i\in I$. (We have ${}^0({}^\ta\btt)={}^\ta\btt$. Indeed, for any 
$f\in\fh_\br$ we have ${}^0({}^0f)={}^0f$.) By symmetry we have also 
$\xi_{\La_i}({}^\ta\btt)\ge\xi_{\La_i}({}^\ta\btt')$ for any $i\in I$. Hence 
$\xi_{\La_i}({}^\ta\btt)=\xi_{\La_i}({}^\ta\btt')$ for any $i\in I$. Thus, the
annihilator of ${}^\ta\btt-{}^\ta\btt'$ in $\fh^*$ contains any $\xi_{\La_i}$ 
with $i\in I$. Since these elements span $\fh^*$, it follows that 
${}^\ta\btt={}^\ta\btt'$. Then for any $\ga\in\cx$ we have 
$\ta(\ga(\btt))=\ta(\ga(\btt'))$. Since this holds for any $\ta$, we deduce 
that $\ga(\btt)=\ga(\btt')$ for any $\ga\in\cx$. Since $\cx$ generates $\fh^*$,
it follows that $\btt=\btt'$.

From 3.26 we see that $K=\{i\in I;\ta(\al_i(\btt))=0\}$. Similarly,
$K'=\{i\in I;\ta(\al_i(\btt'))=0\}$. Since $\btt=\btt'$, it follows that 
$K=K'$. Thus, there exists $g_1\in G$ such that 
$Q'=g_1Qg_1\i,Q'{}^1=g_1Q^1g_1\i$. Replacing $(\si',y',h',\ty')$ by a 
quadruple in the same $G$-orbit, we see that we may assume that 
$Q=Q',Q^1=Q'{}^1$. 

We show that $P'_2$ is automatically contained in $Q$. Assume that 
$P'_2\not\sub Q'=Q$. Let $V\in\ci_K$. Applying 3.19(b) (reformulated with the 
aid of 3.22 and 3.25(b)) to $\si',y',h',\ty',P'_2$, we see that 
$\xi_V({}^\ta\et)>\xi_V({}^\ta\btt')$. Applying 3.19(d) (reformulated with the
aid of 3.22 and 3.25(b)) to $\si,y,h,\ty,P_1$, we see that 
$\xi_V({}^\ta\et)=\xi_V({}^\ta\btt)$. This contradicts the previous inequality
since $\btt=\btt'$. Thus, 
$$P'_2\sub Q'=Q.\tag b$$ 
Note that all of $\si,y,h,\ty,\si',y',h',\ty'$ are contained in $\uq^1$. We
show that these elements satisfy the hypotheses of the lemma 3.16 (with $G$
replaced by $Q^1$).

The analogues of $Q,Q'$ (when $G$ is replaced by $Q^1$) are $Q^1,Q^1$. Now in 
$Q^1$ we have an analogue of $\cp$, namely 
$$\cp'=\{R; R=P\cap Q^1,P\in\cp,P\sub Q\}.$$ 
(See 1.3.) Let $R_1=P_1\cap Q^1,R'_2=P'_2\cap Q^1$. Clearly $R_1\in\cp'$; by 
(b), we have $R_1,R'_2\in\cp'$. Now $\fh$ defined in terms of $Q^1,\cp'$ is 
canonically the same as $\fh$ defined in terms of $G,\cp$. From (i) and (b) we
deduce that
\nl
the image of $\si$ in $\un{R_1}/[\un{R_1},\un{R_1}]=\fh$ coincides with the 
image of $\si'$ in $\un{R'_2}/[\un{R'_2},\un{R'_2}]=\fh$.
\nl
Thus, 3.16 is applicable and $(\si,y,h,\ty)$, $(\si',y',h',\ty')$ are conjugate
under an element of $Q^1$. The lemma is proved.

\proclaim{Lemma 3.31} Assume that $G$ is semisimple. The following two 
conditions on $(y,\si,r)$ are equivalent:

(i) For any $P\in\cp^\si_y$ and any $V\in\ci$ we have 
$\xi_V({}^\ta\ps(P))\ge 0$.

(ii) $Q=G$.
\endproclaim
Assume first that $Q=G$. By 3.28(b) we have $\al_i({}^\ta\btt)=0$ for all
$i\in I$. Hence ${}^\ta\btt=0$. By 3.28(c), for any $P\in\cp^\si_y$ and any 
$V\in\ci$ we have $\xi_V({}^\ta\ps(P)-{}^\ta\btt)\ge 0$, hence
$\xi_V({}^\ta\ps(P))\ge 0$. Thus (ii) implies (i).

Next, assume that $Q\ne G$, that is $K\ne I$. Let $i\in I-K$ and let 
$V\in\ci_{I-\{i\}}$. Then $\xi_V=\sum_{j\in I}z_j\al_j$ (in $\fh^*$) where 
$z_j\ge 0$ for all $j$ and $z_i>0$. By 3.28(a),(b) we have 
$\al_j({}^\ta\btt)\le 0$ for all $j\in I$ and $\al_i({}^\ta\btt)<0$. Hence
$$\xi_V({}^\ta\btt)=\sum_{j\in I}z_j\al_j({}^\ta\btt)<0.$$
Now let $P\in\cp^\si_y$ be such that $P\sub Q$. By 3.28(d) we have
$\xi_V({}^\ta\ps(P)-{}^\ta\btt)=0$ hence $\xi_V({}^\ta\ps(P))<0$. Thus (i) 
implies (ii). The lemma is proved.

\proclaim{Lemma 3.32} Assume that $G$ is semisimple. The following four
conditions on $(y,\si,r)$ are equivalent:

(i) If $P'\in\pa$ and $L'$ is a Levi subgroup of $P'$ such that $\si,y$ are 
contained in $\ul'$, then $P'=G$.

(ii) $y$ is a distinguished nilpotent element of $\fg$ and there exists 
$\hat y\in\ug_{-2r}$ such that $[y,\hat y]=r\i\si$. 

(iii)  For any $P\in\cp^\si_y$ and any $V\in\ci,V\ne\bc$, we have 
$\xi_V(\ps(P))=nr$ where $n\in\bn-\{0\}$.

(iv) For any $P\in\cp^\si_y$ and any $V\in\ci,V\ne\bc$, we have 
$\xi_V({}^\ta\ps(P))>0$.
\endproclaim
The fact that (i) implies (ii) is proved in \cite{\LSQ, 2.5}.

We show that (ii) implies (iii). Assume that (ii) holds. Let 
$\fs=\bc y+\bc\si+\bc\hat y$ (a homomorphic image of $\fs\fl_2(\bc)$). To prove
(iii) it is enough to verify the following statement. 

{\it Let $P\in\cp^\si_y$ and let $V\in\ci,V\ne\bc$. Let $D_P$ be the 
$\up$-stable line in $V$. Then $\si$ acts on $D_P$ as multiplication by $rn$ 
where $n\in\bn-\{0\}$.} (Compare \cite{\LSQ, 2.8}.)
\nl
From the representation theory of $\fs\fl_2$, we see that $\si$ acts on $D_P$ 
as $rn$ where $n\in\bn$. (Use that $y$ acts as $0$ on $D_P$.) Assume that $\si$
acts on $D_P$ as $0$. Then $D_P$ must be stable under $\fs$ which acts on it by
$0$. Now let $P'\in\pa$ be such that $\up'$ is the stabilizer of $D_P$ in $\ug$
(we have $P'\ne G$). Then $\fs\sub\up'$ hence $\fs$ is contained in a Levi 
subalgebra of $\up'$. Thus, $y$ is not distinguished in $\ug$, contradicting 
(ii). Thus (ii) implies (iii).

It is clear that (iii) implies (iv).

Assume now that (iv) holds and (i) does not hold. Then we can find $P'\in\pa$, 
$P'\ne G$, and a Levi subgroup $L'$ of $P'$ such that $\si\in\ul',y\in\ul'$. 
Replacing if necessary $P',L'$ by a $G$-conjugate, we may assume that 
$\si,y,h,\ty$ are all contained in $\ul'$. Let
$A_1=\{P\in\cp^\si_y;P\sub P'\}$. By 3.10 (applied to $P'$ instead of $Q$), we 
have $A_1\ne\emp$. We have $P'\in\cp_J,J\ne I$. Since (iv) holds, we see from
3.31 that $Q=Q^1=G$ and from its proof, that ${}^\ta\btt=0$. Let $V\in\ci_{J}$.
Then $V\ne\bc$. Let $P\in A_1$. By 3.20 (reformulated with the aid of 3.22 and
3.25(b)) we have $\xi_V({}^\ta\ps(P)-{}^\ta\btt)=0$. Since ${}^\ta\btt=0$, it 
follows that $\xi_V({}^\ta\ps(P))=0$. This contradicts (iv) since $V\ne\bc$. 
Thus, (iv) implies (i). The lemma is proved.

\subhead 3.33\endsubhead
Let $E$ be an $\bs$-module of finite dimension over $\bc$. There is a canonical
direct sum decomposition $E=\opl_{\et\in\fh}{}_\et E$ where ${}_\et E$ (a 
{\it weight space}) is the set of all $x\in E$ such that for any $\xi\in\fh^*$,
$\xi:E@>>>E$ is given on ${}_\et E$ by multiplication by $\xi(\et)$ plus a 
nilpotent endomorphism of ${}_\et E$.

\subhead 3.34\endsubhead
The torus $G'=\lan(\si,r)\ran$ in $G\tim\bc^*$ is well defined (see 1.1). The 
fixed point set of the $G'$-action on $\dfg$ (restriction of the $G\tim\bc^*$ 
action) is
$$\dfg^{G'}=\{(y',P)\in\dfg;[\si,y']=2ry',\si\in\up\}.$$ 
Let $pt$ denote a point $(y',P)$ of $\dfg^{G'}$. Since $\si\in\up$, we have 
$\un G'\sub\up\opl\bc$.

Let $k$ be the composition
$$\bs=H^*_{G\tim\bc^*}(\dfg,\bc)@>>>H^*_{G'}(\dfg,\bc)@>>>H^*_{G'}(pt,\bc)
@>>>\bc_{\si,r}$$
where the first map is as in \cite{\LI, 1.4(g)}, the second map is $j^*$ (see
\cite{\LI, 1.4(a)} attached to the imbedding $j:pt@>>>\dfg$ and the third map
is the quotient defined as in 1.13. For $\xi\in\fh^*$ (a subset of $\bs$) we 
have from the definitions:
$$k(\xi)=\xi'(\si,r)\tag a$$
where $\xi'$ is the linear form on $\un G'$ given by the  composition 
$\un G'@>>>\up@>>>\up/[\up,\up]=\fh@>\xi>>\bc$ (the unspecified maps are the
obvious ones).

The image of $\ps:\cp^\si@>>>\fh$ (as in 3.22) is a finite subset $D$ of $\fh$.
From the definitions we have $\xi'(\si,r)=\xi(\ps(P))$ hence (a) can be 
rewritten as follows:
$$k(\xi)=\xi(\ps(P)).\tag b$$
Now let $\tx$ be a subvariety of $\dfg^{G'}$. By 1.11, 
$$H^{G'}_*(\tx,\dcl)=H^*_{G'}\ot H_*(\tx,\dcl)$$
(see \cite{\LII, 1.21}) is naturally an $\bs$-module. Hence
$$A=\bc_{\si,r}\ot_{H^*_{G'}}H_*^{G'}(\tx,\dcl)=H_*(\tx,\dcl)$$
(where $\bc_{\si,r}=H^*_{G'}/\ci^{G'}_{\si,r}$ is as in 1.13) is again an
$\bs$-module.

We have a morphism $\tx@>>>D\sub\fh$ given by $(y',P)\mto\ps(P)$. Consider the
partition $\tx=\sqcup_{\de\in\De}\tx^\de$ of $\tx$ into connected components 
($\De$ is the set of irreducible components of $\tx$.) Now each connected 
component $\tx^\de$ is mapped by $\tx@>>>D$ to a single point of $D$ denoted 
$\ps(\de)$. Since $\tx^\de$ is open and closed in $\tx$ and is $G'$-stable, we
may identify canonically $A=\opl_{\de\in\De}A^\de$ where 
$A^\de=H_*(\tx^\de,\dcl)$. Clearly this direct sum decomposition is compatible
with the $\bs$-module structure.

\proclaim{Lemma 3.35} For any $\et\in\fh$ we have
${}_\et A=\opl_{\de\in\De;\ps(\de)=\et}A^\de$.
\endproclaim
We may assume that $\tx$ is connected. Let $d\in D$ be defined by $d=\ps(P)$ 
for any $(y,P)\in\tx$. Let $\xi\in\fh^*$. We must show that $\xi-\xi(d)$ acts 
nilpotently on $A$. Let $\ti\xi$ be the image of $\xi$ under the composition
$$\align&\bs=H^*_{G\tim\bc^*}(\dfg,\bc)@>>>H^*_{G'}(\dfg,\bc)
@>>>H^*_{G'}(\tx,\bc)\\&=
H^*_{G'}\ot H^*(\tx,\bc)@>>>\bc_{\si,r}\ot H^*(\tx,\bc)\endalign$$
where the first map is as in \cite{\LI, 1.4(g)}, the second map is $\tm^*$ (see
\cite{\LI, 1.4(a)} attached to the imbedding $\tm:\tx@>>>\dfg$, and the third 
map is induced by the quotient defined as in 1.13. The action of $\xi$ on $A$ 
is multiplication by $\ti\xi\in H^*(\tx,\bc)$ on $H_*(\tx,\dcl)$. We have 
$\ti\xi=\ti\xi_0+\ti\xi_>$ where $\ti\xi_0\in H^0(\tx,\bc)$ and
$\ti\xi_>\in\opl_{n>0}H^n(\tx,\bc)$. Clearly, multiplication by $\ti\xi_>$ on 
$H_*(\tx,\dcl)$ is nilpotent. Since $\tx$ is connected, we have 
$\ti\xi_0=c\bold 1$ where $c\in\bc$ and $\bold 1\in H^0(\tx,\bc)$ is the unit 
element of the algebra $H^*(\tx,\bc)$ (which acts as the identity on 
$H_*(\tx,\dcl)$). It is then enough to show that $c=\xi(d)$. Let $pt$ denote a
point $(y',P)$ of $\tx$. Let $j':pt@>>>\tx$ be the imbedding. From the 
definitions, $j'{}^*:H^*(\tx,\bc)@>>>H^*(pt,\bc)=\bc$ carries $\ti\xi$ to 
$k(\xi)$ (as in 4.2), that is to $\xi(\ps(P))=\xi(d)$ (see 4.2(b)). It 
automatically carries $\ti\xi_>$ to $0$ hence it carries $\ti\xi_0$ to 
$\xi(d)$. It also preserves unit elements. Hence it carries $c\bold 1$ to $c$.
Since $\ti\xi_0=c\bold 1$, it follows that $c=\xi(d)$. The lemma is proved.

\subhead 3.36\endsubhead
Let $G'=\lan(\si,r)\ran\sub G\tim\bc^*$. Since $(\si,r)\in\un{M^0(y)}$, we have
$G'\sub M^0(y)$. Let $\cm=E_{y,\si,r}$ be as in 1.13 (recall that the choice of
$G'$ in 1.13 is immaterial; in particular we may take $G'=\lan(\si,r)\ran$). 
For $\rh\in\Irr_0\bam(y,\si)$ we set $\cm_\rh=E_{y,\si,r,\rh}$.

The fixed point set of the $G'$-action on $\cp_y$ (restriction of the 
$M(y)$-action) is just $\cp^\si_y$. By the localization theorem
\cite{\LII, 4.4(b)} (which is applicable in view of the odd vanishing theorem
\cite{\LI, 8.6}), the imbedding $j:\cp^\si_y@>>>\cp_y$ induces an isomorphism 
$$j_!:\bc_{\si,r}\ot_{H^*_{G'}}H_*^{G'}(\cp^\si_y,\dcl)@>\sim>>
\bc_{\si,r}\ot_{H^*_{G'}}H_*^{G'}(\cp_y,\dcl)$$
or equivalently
$$A@>\sim>>\cm.\tag a$$
Here $A$ (as in 3.34) is an $\bs$-module and $\cm$ is an $\bh$-module (in
particular, an $\bs$-module via the obvious algebra homomorphism $\bs@>>>\bh$).
The isomorphism (a) is compatible with $\bs$-module structures. Now the direct
sum decomposition $A=\opl_{\de\in\De}A^\de$ in 3.34 (where ($\De$ is the set of
irreducible components of $\cp^\si_y$) correponds under (a) to a direct sum 
decomposition 
$$\cm=\opl_{\de\in\De}\cm^\de\tag b$$
and we can reformulate 3.35 as follows:

(c) {\it For any $\et\in\fh$ we have}
${}_\et\cm=\opl_{\de\in\De;\ps(\de)=\et}\cm^\de$.
\nl
Here ${}_\et\cm$ are the weight spaces of $\cm$. 

Since the $\bam(y,s)$-action on $\cm$ commutes with the $\bh$-module structure,
each weight space ${}_\et\cm$ is $\bam(y,\si)$-stable. It follows that for
$\rh\in\Irr\bam(y,\si)$, we have

(d) ${}_\et(\cm_\rh)=\Hom_{\bam(y,\si)}(\rh,{}_\et\cm)$.
\nl
Now the $M(y)$-action on $\cp_y$ restricts to an $M(y,\si)$-action on 
$\cp_y^\si$. Hence $\bam(y,\si)=M(y,\si)/M(y,\si)^0$ acts naturally on $\De$.
Let $\baD$ be the set of $\bam(y,\si)$-orbits on $\De$ and let 
$\de\mto\bad$ be the canonical map $\De@>>>\baD$. Let $\rh\in\Irr\bam(y,\si)$.
We have
$$\cm_\rh=\opl_{\ep\in\baD}\cm_\rh^\ep$$ where 
$\cm_\rh^\ep=\Hom_{\bam(y,\si)}(\rh,\opl_{\de\in\De;\bad=\ep}\cm^\de)$. Now
$\ps:\De@>>>\fh$ is constant on $\bam(y,\si)$-orbits hence it induces a map 
$\bar\ps:\baD@>>>\fh$. Then, from (c) and (d) we deduce:

(e) {\it for any $\et\in\fh$ we have}
${}_\et(\cm_\rh)=\opl_{\ep\in\baD;\bar\ps(\ep)=\et}\cm_\rh^\ep$.
\nl
Let $\cp^*=\{P\in\cp;P\sub Q\}$. Let $\De^1$ (resp. $\De^2$) be the set of 
irreducible components of $\cp^\si_y$ that are contained in $\cp^*$ (resp. in 
$\cp-\cp^*$). By 3.21, we have a partition $\De=\De^1\sqcup\De^2$. Let
$$\cm^1=\opl_{\de;\de\in\De^1}\cm^\de,\cm^2=\opl_{\de;\de\in\De^2}\cm^\de.$$
Then $\cm=\cm^1\opl\cm^2$. Now the summands $\cm^1,\cm^2$ of $\cm$ are stable 
under $\bam(y,\si)$ (since $M(y,\si)\sub Q\tim\bc^*$, see 3.15). Hence, setting

$\cm_\rh^1=\Hom_{\bam(y,\si)}(\rh,\cm^1), 
\cm_\rh^2=\Hom_{\bam(y,\si)}(\rh,\cm^2)$, 
\nl
we have $\cm_\rh=\cm^1_\rh\opl\cm^2_\rh$.

\proclaim{Lemma 3.37} (a) Each of $\cm^1$ and $\cm^2$ is a sum of weight spaces
of $\cm$. 

(b)  Each of $\cm^1_\rh$ and $\cm^2_\rh$ is a sum of weight spaces of 
$\cm_\rh$. 
\endproclaim
Using the inclusion $\cm^\de\sub{}_{\ps(\de)}\cm$ (see 3.36(c)) we see that

(c) $\cm^1\sub\cm^{(1)},\cm^2\sub\cm^{(2)}$ where

$\cm^{(1)}=\sum_{\de;\de\in\De^1}{}_{\ps(\de)}\cm,
\cm^{(2)}=\sum_{\de;\de\in\De^2}{}_{\ps(\de)}\cm$.
\nl
We show that 

(d) $\cm^{(1)}\cap\cm^{(2)}=0$. 
\nl
Since the weight spaces of $\cm$ form a direct sum, it is enough to show that 
$\ps(\de)\ne\ps(\de')$ for any $\de,\de'\in\De$ such that 
$\de\sub\cp^*,\de'\sub\cp-\cp^*$ or equivalently, that $\ps(P)\ne\ps(P')$ for
any $P,P'\in\cp^\si_y$ such that $P\sub Q,P'\not\sub Q$. 

Let $V\in\ci_K$. By 3.19(b),(d) we have $\ta(\nu_V(P))/\ta(r)=b_V$, 
$\ta(\nu_V(P'))/\ta(r)>b_V$, hence $\nu_V(P)\ne\nu_V(P')$. Hence 
$\xi_V(\ps(P))\ne\xi_V(\ps(P'))$. Hence $\ps(P)\ne\ps(P')$, as desired, and (d)
is proved.

Since $\cm=\cm^1\opl\cm^2$, we see from (c),(d) that $\cm^1=\cm^{(1)}$, 
$\cm^2=\cm^{(2)}$. Since each $\cm^{(1)},\cm^{(2)}$ is a sum of weight spaces 
of $\cm$, (a) follows.

We prove (b). Let $k\in\{1,2\}$. By (a) we have 
$\cm^k=\opl_{n=1}^N({}_{\et_n}\cm)$
 where $\et_n$ are distinct elements of $\fh$.
It follows that
$$\cm^k_\rh=\Hom_{\bam(y,\si)}(\rh,\cm^k)=\opl_{n=1}^N
\Hom_{\bam(y,\si)}(\rh,{}_{\et_n}\cm)=\opl_{n=1}^N({}_{\et_n}\cm_\rh).$$
The lemma is proved.

\subhead 3.38\endsubhead
Replacing $G,y,\si,r$ by $Q^1,y,\si,r$ in the definition of 
$W,\fh,\bs,\bh,\cp_y,\cp^\si_y,\dcl,\cm$ we obtain (as in 1.16)
$W_K,\fh,\bs,\bh',\cp'_y,\cp'{}^\si_y,\dcl,\cm'$. Since $M(y,\si)\sub Q$, the
analogue of $M(y,\si)$ for $Q^1$ instead of $G$ is a quotient of $M(y,\si)$ by
a unipotent normal subgroup. Hence $\bam(y,\si)$ defined in terms of $Q^1$ is 
the same as that defined in terms of $G$. Hence for $\rh\in\Irr\bam(\si,y)$, we
can define $\cm'_\rh$ in terms of $\cm'$ in the same way as $\cm_\rh$ was 
defined in terms of $\cm$. By 1.18 we have an isomorphism
$$\Psi:\bh\ot_{\bh'}\cm'@>\sim>>\cm.\tag a$$
Then $x\mto\Psi(1\ot x)$ is a map $\cm'@>>>\cm$.

(b) This is an isomorphism of $\cm'$ onto $\cm^1$. 
\nl
Indeed, using the localization theorem \cite{\LII, 4.4(b)}, $\cm'@>>>\cm$ may
be identified with the map

$\bc_{\si,r}\ot_{H^*_{G'}}H_*^{G'}(\cp'{}^\si_y,\dcl)@>>>
\bc_{\si,r}\ot_{H^*_{G'}}H_*^{G'}(\cp^\si_y,\dcl)$
\nl
induced by the obvious inclusion $\cp'{}^\si_y\sub\cp^\si_y$ whose image is the
open and closed subset $\cp^\si_y\cap\cp^*$ of $\cp^\si_y$.

Now for any $\rh\in\Irr\bam(\si,y)$, (a) induces an isomorphism
$$\bh\ot_{\bh'}\cm'_\rh@>\sim>>\cm_\rh\tag c$$
and (b) induces an isomorphism
$$\cm'_\rh@>\sim>>\cm^1_\rh.\tag d$$
Since $\bh'@>>>\bh$ is injective and $\bh$ is a free $\bh'$-module, from (c) we
deduce that 
$$\cm'_\rh\ne 0\lra\cm_\rh\ne 0.\tag e$$

\subhead 3.39. Proof of Theorem 1.15(a)\endsubhead 
Note that in the setup of 1.15(a) we have $E_{y,\si,r}\ne 0$ hence we are also
in the setup of 3.9.

By 3.38(e) we have $\cm'_\rh\ne 0$. As in the proof of 3.16 we see that $y$ 
belongs to the unique open orbit of $Z(\si)\cap Q^1$ in $\fg_{2r}\cap\fq^1$. 
Hence we may apply \cite{\LI, 8.17(b)} to $Q^1$ instead of $G$ and conclude 
that $\cm'_\rh$ is a simple $\bh'$-module.

Let $F$ be a proper $\bh$-submodule of $\cm_\rh$. If $F\cap\cm^1_\rh=\cm^1_\rh$
then, since $\cm^1_\rh$ generates $\cm_\rh$ as a $\bh$-module (see 3.28(c),(d))
it would follow that $F$ generates $\cm_\rh$ as a $\bh$-module hence 
$F=\cm_\rh$ contradicting the assumption that $F$ is proper.

Thus we must have $F\cap\cm^1_\rh\ne\cm^1_\rh$. Since $\cm^1_\rh$ is a simple
$\bh'$-module and $F\cap\cm^1_\rh$ is a proper $\bh'$-submodule of $\cm^1_\rh$,
it follows that 

(a) $F\cap\cm^1_\rh=0$. 
\nl
For any $\et\in\fh$ we have ${}_\et F\sub {}_\et(\cm_\rh)$ and by 3.37(b),
${}_\et(\cm_\rh)$ is contained either in $\cm^1_\rh$ or in $\cm^2_\rh$. Thus,
we have either ${}_\et F\sub\cm^1_\rh$ or ${}_\et F\sub\cm^2_\rh$. The first 
alternative cannot occur if ${}_\et F\ne 0$ by (a). Thus, we have
${}_\et F\sub\cm^2_\rh$ for any $\et$ hence $F\sub\cm^2_\rh$. It follows that 
the sum of all proper $\bh$-submodules of $\cm_\rh$ is contained in $\cm^2_\rh$
and thus it is itself a proper submodule (since $\cm^1_\rh$ is $\ne 0$, being 
isomorphic to $\cm'_\rh$). This proves 1.15(a).

\subhead 3.40\endsubhead
Let $\rh\in\Irr_0\bam(\si,y)$. Let $\cm_{\rh,max}$ be the unique maximal 
$\bh$-submodule of $\cm_\rh$. Recall that $\cm_{\rh,max}\sub\cm^2_\rh$. It
follows that

(a) {\it The obvious map $\cm'_\rh=\cm^1_\rh@>>>\cm_\rh/\cm_{\rh,max}$ is 
injective.}
\nl
Let $V\in\ci_K$. Let 
$$X=\{\et\in\fh;{}_\et(\cm_\rh/\cm_{\rh,max})\ne 0, 
\ta(\xi_V(\et))/\ta(r) \text{ is minimum possible }\}.$$
Note that $X$ is well defined (see the proof of 3.37) and the minimum value in
the definition of $X$ is $b_V$. From the proof of 3.37 we see also that

(b) {\it The image of the map in (a) equals}
$\sum_{\et\in X}{}_\et(\cm_\rh/\cm_{\rh,max})$.

\subhead 3.41. Proof of injectivity in Theorem 1.15(b)\endsubhead
Let $y^!,\si^!,h^!,\ty^!$ be another quadruple like $y,\si,h,\ty$ (with the
same $r$). Let $\rh\in\Irr_0\bam(\si,y)$, $\rh^!\in\Irr_0\bam(\si^!,y^!)$.
Define 
$$Q^!,Q^1{}^!,\cm^!_{\rh^!},\cm'{}^!_{\rh^!},\cm^1{}^!_{\rh^!},
\cm^!_{\rh^!,max},\ps'$$ 
in terms of $y^!,\si^!,h^!,\ty^!,\rh^!$ in the same way as 
$$Q,Q{}^!,\cm_\rh,\cm'_\rh,\cm^1_\rh,\cm_{\rh,max},\ps$$
were defined in terms of $y,\si,h,\ty,\rh$. Assume that 
$$\cm_\rh/\cm_{\rh,max}\cong\cm^!_{\rh^!}/\cm^!_{\rh^!,max}\tag a$$ 
as $\bh$-modules. We show that there exists $g\in G$ which conjugates 
$y,\si,h,\ty,\rh$ to $y^!,\si^!,h^!,\ty^!,\rh^!$.

We can find $\et\in\fh$ such that ${}_\et(\cm'_\rh)\ne 0$. By 3.40(a), we have
${}_\et(\cm_\rh/\cm_{\rh,max})\ne 0$ hence by (a),
${}_\et(\cm^!_{\rh^!}/\cm^!_{\rh^!,max})\ne 0$.

It follows that ${}_\et\cm^1\ne 0$, ${}_\et(\cm^!_{\rh^!})\ne 0$.

Using 3.36(c) we deduce that there exist $P_1\in\cp^\si_y$ and
$P'_2\in\cp^{\si'}_{y'}$ such that $P_1\sub Q$ and

$\ps(P_1)=\et=\ps^!(P'_2)$
\nl 
By symmetry, there exist $P'_1\in\cp^{\si'}_{y'}$, $P_2\in\cp^\si_y$ and 
$\et'\in\fh$ such that $P'_1\sub Q$ and

$\ps^!(P'_1)=\et'=\ps(P_2)$.
\nl
Then the assumptions of 3.30 are verified and we see that there exists $g\in G$
which conjugates $y,\si,h,\ty$ to $y^!,\si^!,h^!,\ty^!$. Thus, we may assume 
that $y^!=y,\si^!=\si,h^!=h,\ty^!=\ty$. Then 
$$\align&Q^!=Q,Q^1{}^!=Q^1,\cm^!_{\rh^!}=\cm_{\rh^!},
\cm'{}^!_{\rh^!}=\cm'_{\rh^!},\\&
\cm^1{}^!_{\rh^!}=\cm^1_{\rh^!},\cm^!_{\rh^!,max}=\cm_{\rh^!,max},\ps'=\ps.
\endalign$$
Let $V\in\ci_K$. Consider an $\bh$-linear isomorphism 
$\cm_\rh/\cm_{\rh,max}@>\sim>>\cm_{\rh^!}/\cm_{\rh^!,max}$. This clearly 
carries the subspace

$\sum_{\et\in X}{}_\et(\cm_\rh/\cm_{\rh,max})$
\nl
($X$ as in 3.40) onto the subspace

$\sum_{\et\in X}{}_\et(\cm_{\rh^!}/\cm_{\rh^!,max})$
\nl
(here $X$ is as in 3.40, and the analogous set for $\cm_{\rh^!}$ is again $X$)
hence, by 3.40, it carries the subspace $\cm^1_\rh$ isomorphically onto the 
subspace $\cm^1_{\rh^!}$. Hence it induces an isomorphism of $\bh'$-modules
$\cm'_\rh@>\sim>>\cm^1_{\rh^!}$. As in the proof of 3.16 we see that $y$ 
belongs to the unique open orbit of $Z(\si)\cap Q^1$ in $\fg_{2r}\cap\fq^1$. We
now apply \cite{\LI, 8.17(c)} to $Q^1$ instead of $G$ and deduce that 
$\rh^!=\rh$. This completes the proof.

\subhead 3.42. Proof of surjectivity in Theorem 1.15(b)\endsubhead
A statement close to the surjectivity in Theorem 1.15(b) was stated in 
\cite{\LI, 8.15} but the proof given there has an error in line -7 of p.199 
("Since $H^*_{M^0(y)}/I''$ is an artinian $\bc$-algebra..."). (I thank David 
Vogan for pointing out that error). In the part of the proof preceding that 
line it was shown that:

(a) if $\cn$ is a simple $\bh$-module, then there exists $y\in\fg_N$, an ideal
$J$ of finite codimension in $\ca=H^*_{M^0(y)}$ and a non-zero $\bh$-linear map
$X/JX@>>>\cn$ where $X=H_*^{M^0(y)}(\cp_y,\dcl)$. 
\nl
We continue the proof starting from (a). We have $\dim_\bc(X/JX)<\infty$ (this
can be seen from \cite{\LI, 7.2} using \cite{\LI, 8.6}). Let $X/JX@>>>X'$ be 
the largest semisimple quotient of the $\bh$-module $X/JX$. Then $X'$ inherits
from $X/JX$ an $\ca$-module structure. Clearly, there exists a non-zero 
$\bh$-linear map $X'@>>>\cn$. Let $X'_{\cn}$ be the $\cn$-isotypical part of 
the $\bh$-module $X'$. Then $X'_{\cn}$ is an $\ca$-submodule and 
$X'_{\cn}\ne 0$. Let $rad J=\{x\in\ca;x^n\in J \text{ for some } n\ge 1\}$. The
elements of $rad(J)$ act on $X'_{\cn}$ as commuting nilpotent elements. Hence
$rad(J)X'_{\cn}\ne X'_{\cn}$ (since $X'_{\cn}\ne 0$). Hence there exists a
non-zero $\bh$-linear map $X'_{\cn}/rad(J)X'_{\cn}@>>>\cn$. Now 
$X'_{\cn}/rad(J)X'_{\cn}$ is a direct summand of the $\bh$-module $X'/rad(J)X'$
hence there exists a non-zero $\bh$-linear map $X'/rad(J)X'@>>>\cn$. The 
canonical map $X/rad(J)X@>>>X'/rad(J)X'$ is surjective hence there exists a 
non-zero $\bh$-linear map $X/rad(J)X@>>>\cn$. Since $rad(rad(J))=rad(J)$ we see
that we may assume that $J=rad(J)$. Then the commutative algebra $A/JA$ is a 
finite direct sum of copies of $\bc$. Let $I_1,\dots,I_k$ be the maximal ideals
of $A$ that contain $J$. We have $A/J@>\sim>>A/I_1\opl\dots\opl A/I_k$ and
$X/JX=X/I_1X\opl\dots X/I_kX$ (as $\bh$-modules). Hence there exists 
$j\in[1,k]$ and a non-zero $\bh$-linear map $X/I_jX@>>>\cn$. Hence we may 
assume that $J$ is a maximal ideal of $\ca$. Then there exists a semisimple 
element $(\si,r)\in\un{M^0(y)}$ such that $J=\cj^{M^0(y)}_{\si,r}$ (see 1.13).
Then $X/JX=E_{y,\si,r}$ and we see that there exists a non-zero $\bh$-linear 
map $E_{y,\si,r}@>>>\cn$. Now $E_{y,\si,r}=\opl_\rh\rh\ot E_{y,\si,r,\rh}$
where $\rh$ runs over $\Irr_0\bam(y,\si)$. Hence there exists
$\rh\in\Irr_0\bam(y,\si)$ and a non-zero $\bh$-linear map 
$E_{y,\si,r,\rh}@>>>\cn$. This map is surjective since the $\bh$-module $\cn$ 
is simple. It follows that the induced map $\bae_{y,\si,r,\rh}@>>>\cn$ is an 
isomorphism. This completes the proof of Theorem 1.15.

\subhead 3.43. Proof of Theorem 1.21\endsubhead
Assume that 1.21(iii) holds. Replacing $h,\ty$ by $\ph(h_0),\ph(f_0)$, we see 
that $G=Q^1=Q$. Using 3.31 we see that 3.31(i) holds. Using 3.36(c) we deduce 
that for any $\et\in\fh$ such that ${}_\et\cm\ne 0$ and any $V\in\ci$ we have 
$\ta(\xi_V(\et))/\ta(r)\ge 0$. Hence $\cm$ is $\ta$-tempered. Hence $\cm_\rh$
is $\ta$-tempered and 1.21(i) holds.

Clearly, if 1.21(i) holds, then 1.21(ii) holds.

Assume now that 1.21(iii) does not hold. Using 3.31 we see that $Q\ne G$ that
is $i\in I-K$. Let $V\in\ci_{I-\{i\}}$. Let $\et\in\fh$ be such that 
${}_\et\cm'_\rh\ne 0$. Then, by 3.40(a), we have 
${}_\et(\cm_\rh/\cm_{\rh,max})\ne 0$. We can find $P\in\cp^\si_y$ such that 
$P\sub Q$ and $\et=\ps(P)$. By the second paragraph in the proof of 3.31 we 
have $\ta(\xi_V(\ps(P)))/\ta(r)<0$ hence $\ta(\xi_V(\et))/\ta(r)<0$. It follows
that $\cm_\rh/\cm_{\rh,max}$ is not $\ta$-tempered. Thus 1.21(ii) does not 
hold. 

Thus the three conditions in 1.21 are equivalent. If these conditions are
satisfied then, as we have seen above, we have $G=Q$. As in the proof of 3.16 
we see that $y$ belongs to the unique open orbit of $Z(\si)$ in $\fg_{2r}$. 
Hence we may apply \cite{\LI, 8.17(b)} and conclude that $\cm_\rh$ is a simple
$\bh$-module. This completes the proof of 1.21.

\subhead 3.44. Proof of Theorem 1.22\endsubhead
The equivalence of 1.22(i) and 1.22(ii) follows from 3.32.

Assume that 1.22(ii) holds. Using 3.32 we see that 3.32(iii) holds. Using 
3.36(c) we deduce that for any $\et\in\fh$ such that ${}_\et\cm\ne 0$ and any 
$V\in\ci,V\ne\bc$ we have $\xi_V(\et)=nr$ where $n\in\{1,2,3,\dots\}$. Hence 
1.22(v) holds.

Clearly, if 1.22(v) holds, then 1.22(iv) holds.

By 1.21, conditions 1.22(iii) and 1.22(iv) are equivalent.

Assume that 1.22(iv) holds and that 1.22(i) does not hold. Then we can find 
$P'\in\pa$, $P'\ne G$ and a Levi subgroup $L'$ of $P'$ such that
$y\in\ul',\si\in\ul'$. We may assume that $h\in\ul',\ty\in\ul'$. By 3.10, we 
have $\{P\in\cp^\si_y;P\sub P'\}\ne\emp$. In particular, $P'\in\cp_J$ for some
$J\sub I$, $J\ne I$. In particular $L'$ inherits a natural cuspidal datum
(as in 1.3) and in terms of this we can define $\ti{\bh},\ti{\cm}$ in the same
way as $\bh,\cm$ were defined in terms of the cuspidal datum of $G$. We have a
natural imbedding 
$$j_0:\ti{\cm}@>>>\cm.$$
Let $\fn'=\un{U_{P'}}$ and let ${}_y\fn'=\cok(\ad(y):\fn'@>>>\fn')$. We show
that

(a) $\ad(\si)-2r:{}_y\fn'@>>>{}_y\fn'$ is invertible. 
\nl
By 1.21 (and its proof) we have $Q=Q^1=G$. Thus, 
${}_n\fg_a\ne 0\impl\ta(a)/\ta(r)=n$. We must show that any eigenvalue of 
$\ad(\si)-2r$ on $\cok(\ad(y):\fn'@>>>\fn')$ is $\ne 0$. Now the last cokernel
is a quotient of $\opl_{n,a;n\le 0}({}_n\fg_a)$. Hence it suffices to show that
if ${}_n\fg_a\ne 0$ and $n\le 0$ then $a-2r\ne 0$. But 

$\ta(a-2r)/\ta(r)=(n\ta(r)-2\ta(r))/\ta(r)=n-2\le-2$
\nl
hence $a-2r\ne 0$, as required. This proves (a).

We see that 1.18 is applicable (with $P',L'$ instead of $Q,Q^1$). We deduce
that

(b) $j_0(\ti{\cm})$ generates the $\bh$-module $\cm$.
\nl
We can find a surjective $\bh$-linear map $p:\cm@>>>\cm_\rh$. The composition 
$\ti{\cm}@>j_0>>\cm@>p>>\cm_\rh$ is non-zero. (Otherwise, $j_0(\ti{\cm})$ would
be contained in the proper $\bh$-submodule $\Ker(p)$ of $\cm$ contradicting 
(b).) Since $pj_0$ is $\bs$-linear, it follows that there exists $\et\in\fh$ 
such that ${}_\et\ti{\cm}\ne 0$ and ${}_\et\cm_\rh\ne 0$. Hence there exists 
$P\in\cp^\si_y$ such that $P\sub P'$ and $\ps(P)=\et$. By the argument in the 
last paragraph of the proof of 3.32, we have $\ta(\xi_V(\ps(P)))/\ta(r)=0$, 
that is, $\ta(\xi_V(\et))/\ta(r)=0$ where $V\in\ci_J$. This contradicts 
1.22(iv) since $V\ne\bc$. We have proved that if 1.22(iv) holds then 1.22(i) 
holds. 

This completes the proof of 1.22.

\enddocument